\documentclass[pdflatex,sn-mathphys-num]{sn-jnl}

\usepackage{graphicx}
\usepackage{amsmath,amssymb,amsfonts}
\usepackage{amsthm}
\usepackage{xcolor}
\usepackage{float}
\usepackage{caption}
\usepackage{tikz}
\usepackage{pgfplots}
\pgfplotsset{compat=1.18}

\newtheorem{thm}{Theorem}[section]
\newtheorem{lem}[thm]{Lemma}
\newtheorem{cor}[thm]{Corollary}
\newtheorem{rem}[thm]{Remark}
\newtheorem{pro}[thm]{Proposition}
\newtheorem{defi}[thm]{Definition}
\newtheorem{example}[thm]{Example}

\begin{document}
	
	\title[Self-similar solutions to the inverse mean curvature flow in the Lorentz--Minkowski plane]{Self-similar solutions to the inverse mean curvature flow in the Lorentz--Minkowski plane}
	
	\author{\fnm{Li} \sur{Ou}}\email{ouli@swufe.edu.cn}

	\affil{\orgdiv{School of Mathematics}, \orgname{Southwestern University of Finance and Economics}, \orgaddress{\street{No.555, Liutai Road, Wenjiang}, \city{Chengdu}, \postcode{611130}, \state{Sichuan Province}, \country{China}}}
	
	\abstract{We classify the self-similar solutions of the inverse mean curvature flow in the Lorentz--Minkowski plane $\mathbb{R}^{1,1}$. Using hyperbolic numbers, the self-similarity condition becomes a linear ODE system, and all solutions are obtained in closed form with explicit smoothness and completeness criteria. Translating solitons occur in three causal types, and dilation--rotation solitons are governed by the discriminant $\Delta=a^2-4(b-1)$, as in the Euclidean case. When $a^2=b^2$ the translation cannot be eliminated, owing to the zero divisors of the hyperbolic numbers, and a genuinely Lorentzian class of mixed solitons arises, whose complete smooth members form a one-parameter family. The planar classification is then extended to translationally invariant solutions in $\mathbb{R}^{n,1}$, and rotationally invariant solutions are reduced to an explicit ODE system.}
	
	\keywords{inverse mean curvature flow, self-similar solutions, Lorentz--Minkowski plane, hyperbolic numbers, translating solitons, hyperbolic logarithmic spirals}
	
	\pacs[MSC Classification]{53E10, 53A10, 53C42, 34A30}	
	\maketitle	
\section{Introduction}\label{sec1}
The inverse mean curvature flow (IMCF) moves a hypersurface in the direction of its unit normal $\nu$ with speed $1/H$,
\begin{equation}\label{eq:imcf}
	\frac{\partial X}{\partial t} = \frac{1}{H}\nu,
\end{equation}
where $H$ is the mean curvature. Equation~\eqref{eq:imcf} is a quasilinear parabolic equation which degenerates wherever $H$ vanishes, and this degeneracy is the main analytic difficulty of the flow. In contrast to the mean curvature flow, which contracts hypersurfaces, the inverse flow is expanding: convex hypersurfaces in Euclidean space become round after rescaling~\cite{MR1085136}, and the same holds for star-shaped hypersurfaces~\cite{MR1082861,Gerhardt1990}. The flow entered mathematical physics through the proof by Huisken and Ilmanen of the Riemannian Penrose inequality~\cite{Huisken2001,huisken1997riemannian} and the work of Bray and Neves on $3$-manifolds~\cite{MR2081443}.

Self-similar solutions are the solutions which evolve, up to a similarity of the ambient space, without changing their shape. They are the natural candidates for the rescaling limits of the flow at a singularity, so that a classification of them amounts to a list of the possible blow-up models, and the explicit solutions show which of these models are complete and which become singular after a finite parameter time. In Euclidean space the planar theory is well understood: the translating solitons of IMCF in $\mathbb{R}^2$ are cycloids \cite{Drugan2016}, the homothetic solitons include the logarithmic spirals \cite{Drugan2016,MR1793674}, and Chang~\cite{Chang2020} classified all self-similar motions of curves under IMCF through the equation
\begin{equation*}
	c\langle\gamma, T\rangle + d\langle\gamma, N\rangle = -\frac{1}{k},
\end{equation*}
where $c$ and $d$ are the rotation and the expansion rate, obtaining a trichotomy governed by the discriminant $c^{2}-4(1-d)$.

The Lorentzian theory of IMCF is far less developed. The reason is the indefiniteness of the metric: space-like, time-like and light-like curves behave in fundamentally different ways, and even the definition of the curvature depends on the causal type of the curve. Self-similar solutions of the mean curvature flow in the Minkowski plane $\mathbb{R}^{1,1}$ were classified by Halldorsson \cite{Halldorsson2015}, and da Silva and Tenenblat studied the curvature and the inverse curvature flow on the light cone \cite{daSilva2023}, but a classification of the self-similar solutions of IMCF in $\mathbb{R}^{1,1}$ which incorporates all causal types and the full similarity group, hyperbolic rotations and translations included, has been lacking. The present paper provides this classification, and shows that it has direct consequences in every dimension.

\subsection{Setting and main results}

The Minkowski plane $\mathbb{R}^{1,1}$ is $\mathbb{R}^2$ with the inner product $\langle u,v\rangle_m=u_1v_1-u_2v_2$. We identify it with the ring of hyperbolic numbers $z=x+hy$, $h^2=1$, $h\ne\pm1$, so that $\langle z_1,z_2\rangle_m=\Re(\bar z_1z_2)$ and the hyperbolic rotations (Lorentz boosts) become the multiplications by $e^{h\theta}=\cosh\theta+h\sinh\theta$ \cite{hypernumber}. For light-like curves the classical curvature theory degenerates, because the tangent and the normal directions coincide, so we consider curves which are everywhere space-like or everywhere time-like; since the reflection $z\mapsto hz$ interchanges the two types and carries solutions of the flow to solutions, it suffices to treat space-like curves. A self-similar motion of a curve $\gamma$ is a family
\begin{equation*}
	\hat{\gamma}(s,t)=g(t)e^{hf(t)}\gamma(s)+\mathbf{G}(t),
\end{equation*}
combining a dilation $g$, a hyperbolic rotation $f$ and a translation $\mathbf{G}$. We show in Section~\ref{sec3} that $\hat{\gamma}$ solves IMCF if and only if $\gamma$ satisfies
\begin{equation}\label{eq:meqintro}
	a\langle \gamma, \mathbf{T}\rangle_m-b\langle \gamma, \mathbf{N}\rangle_m-\langle C, \mathbf{N}\rangle_m=\frac{1}{k},
\end{equation}
where $a=f'(0)$, $b=g'(0)$ and $C=\mathbf{G}'(0)$ are the rates of rotation, dilation and translation, $\mathbf{T},\mathbf{N}$ are the unit tangent and normal and $k$ is the curvature. In terms of the hyperbolic tangent angle $\theta$, defined by $\mathbf{T}=e^{h\theta}$, equation~\eqref{eq:meqintro} becomes a linear system of ordinary differential equations with constant coefficients, and all its solutions are obtained in closed form. A solution is called \emph{smooth} if $\frac{1}{k}$ never vanishes, and \emph{complete} if its arc length diverges at both ends (Definition~\ref{def:complete}); a solution may be complete and still have cusps.

If $a^2\ne b^2$, the hyperbolic number $b+ha$ is invertible and the translation term in \eqref{eq:meqintro} is removed by a change of origin, exactly as in the Euclidean case. If $a^2=b^2$, however, $b+ha$ is a zero divisor and the translation cannot be eliminated. This is the point where the Lorentzian theory departs from the Euclidean one, and it produces a class of solutions with no Euclidean counterpart. Our main result is the following.

\begin{thm}[Main Theorem]\label{thm:main}
	Up to homotheties and isometries of $\mathbb{R}^{1,1}$, the space-like solutions of \eqref{eq:meqintro} are the following.
	\begin{enumerate}
		\item[\rm(I)] For $a=b=0$: three translating solitons, one for each causal type of $C$ (Theorem~\ref{translating thm}).
		\item[\rm(II)] For $C=0$, which is no restriction when $a^2\ne b^2$: the dilation--rotation solitons, classified by the discriminant $\Delta=a^2-4(b-1)$, with one, two or four solutions according as $\Delta<0$, $\Delta=0$ or $\Delta>0$ (Theorem~\ref{main thm1}).
		\item[\rm(III)] For $a^2=b^2\ne0$ and $C\ne0$: the mixed solitons, with two solutions if $a=b=2$ and four for every other value of $a=b$ (Theorem~\ref{thm:6.1}).
	\end{enumerate}
	All solutions are given in closed form, and for each of them it is decided whether it is smooth and whether it is complete.
\end{thm}

The completeness results are as follows. Among the translators only the time-like one is complete and smooth. In the hyperbolic case $\Delta>0$ the two non-spiral solutions $\gamma_1,\gamma_2$ are complete if and only if $b\le1$ (for $b\ne\pm a$), and exactly one of them is smooth, according to the sign of $b^2-a^2$; in particular every pure hyperbolic rotation and every pure dilation with $b<1$ admits a complete smooth solution. In the elliptic and parabolic cases completeness forces $a=0$, and elliptic solutions always have infinitely many cusps. Among the mixed solitons the complete smooth ones are exactly the curves $\gamma_+$ of \eqref{eq:gpm} with $a\le1$ and the standard hyperbola at $a=1$ (Corollary~\ref{cor:6.3}); since the rate $a$ cannot be normalized away, they form a genuine one-parameter family.

Several features of this classification are specific to the Lorentzian setting. Whereas the cycloid is the unique translating soliton in $\mathbb{R}^2$, translating solitons in $\mathbb{R}^{1,1}$ come in three causal types, and the time-like one is the Wick rotation of the cycloid. The hyperbolic logarithmic spirals $\gamma=(l-h)e^{(l+h)\theta}$ (Proposition~\ref{hyper log spiral}) play the role of the logarithmic spirals, and the discriminant $\Delta$ is the exact counterpart of Chang's $c^2-4(1-d)$; but the completeness properties differ, because the arc length is measured by the indefinite metric. The two lines $b=\pm a$ are exceptional already in the dilation--rotation case: on them equation \eqref{eq:meqintro} acquires a translational symmetry along a light-like direction (Remark~\ref{rem:nulltrans}), and this is the mechanism behind class (III). Finally, for the inverse mean curvature flow the rates $a,b$ are invariant under rescaling, in contrast with the mean curvature flow where $|b|$ can be normalized to $1$ \cite{Halldorsson2015}; this is why class (III) contains a one-parameter family rather than isolated examples.

The planar classification has direct consequences in arbitrary dimension (Section~\ref{sec7}). Writing $\mathbb{R}^{n,1}=\mathbb{R}^{1,1}\times\mathbb{R}^{n-1}$, the cylinder $\gamma\times\mathbb{R}^{n-1}$ over a planar soliton is a self-similar solution of IMCF in $\mathbb{R}^{n,1}$, and conversely every self-similar solution invariant under the translations of $\mathbb{R}^{n-1}$ is such a cylinder. Hence Theorem~\ref{thm:main} classifies the $\mathbb{R}^{n-1}$-invariant self-similar solutions in all dimensions $n\ge2$ (Theorem~\ref{thm:cyl}); the cylinder over the complete translating soliton is the analogue, for IMCF in $\mathbb{R}^{n,1}$, of the grim-reaper cylinder of the Euclidean mean curvature flow. For hypersurfaces invariant under rotations about a time-like axis the mean curvature acquires the extra term $(n-1)y'/r$, the hyperbolic rotations are excluded, and the problem reduces to the explicit system \eqref{eq:rotode}, whose umbilic solution is the hyperbolic space $\mathbb{H}^n(\rho)$ expanding at the rate $\rho(t)=\rho_0e^{t/n}$ (Example~\ref{ex:hyperboloid}).

The paper is organized as follows. Section~\ref{sec2} reviews the hyperbolic numbers and the Minkowski plane, and Section~\ref{sec3} derives equation \eqref{eq:meqintro}. Sections~\ref{sec4}, \ref{sec5} and \ref{sec6} classify the translating, the dilation--rotation and the mixed solitons respectively, with their completeness, smoothness and asymptotic behaviour. Section~\ref{sec7} derives the higher-dimensional consequences.

	\section{The hyperbolic number plane}\label{sec2}
	To study the geometric properties of self-similar solutions in non-Euclidean spaces, we first review the fundamental structure of the Minkowski plane $\mathbb R^{1,1}$ and its connection to hyperbolic numbers.
	
	Let $\mathbb R^{1,1}$ be the 2-dimensional Minkowski plane.  The Minkowski inner product is
	\begin{equation*}
		\langle (x_1,y_1),(x_2,y_2)\rangle_m:=x_1x_2-y_1y_2.
	\end{equation*}
	A nonzero vector $\mathbf u \in \mathbb R^{1,1}$ is called space-like if $\langle \mathbf u,\mathbf u \rangle_m>0$, time-like if $\langle \mathbf u,\mathbf u \rangle_m<0$, and light-like (or null) if $\langle \mathbf u,\mathbf u \rangle_m=0$. Two vectors $\mathbf u_1, \mathbf u_2$ are said to be orthogonal if and only if $\langle \mathbf u_1,\mathbf u_2 \rangle_m=0,$ and $||\mathbf u||_m=\sqrt{|\langle\mathbf u,\mathbf u\rangle_m|}$ is called the \text{Minkowski norm}.
	
	For ease of computation in the Minkowski plane, we utilize hyperbolic numbers (also called split-complex numbers) throughout this paper. Each point $(x,y)\in \mathbb R^{1,1}$ corresponds to a hyperbolic number $z=x+hy$, where the hyperbolic unit  $h$ satisfies  $h^2=1$ but $h\ne \pm 1$ as opposed to $i^2=-1$ for the complex numbers. Thus, the hyperbolic numbers $\mathbb D=\mathbb R[h]/(h^{2}-1)$ form an extension of the real numbers by adjoining the unipotent element $h$, analogous to how the complex numbers $\mathbb C=\mathbb R[i]$  extend the reals by including the imaginary unit $i$. The addition and multiplication of hyperbolic numbers are defined as follows:		
	
	$$(x_1+hy_1)+(x_2+hy_2)=(x_1+x_2)+h(y_1+y_2),$$
	$$(x_1+hy_1)\cdot (x_2+hy_2)=(x_1 x_2+y_1 y_2)+h(x_1y_2+x_2y_1).$$
	Just as complex numbers $a+ib$ can be represented as the $2\times 2$ matrices
	$\left(\begin{array}{cc}
		a & -b \\
		b & a
	\end{array}\right),$ hyperbolic numbers $x+hy$ correspond to symmetric $2\times 2$ matrices
	$\left(\begin{array}{cc}
		x & y\\
		y & x
	\end{array}\right).
	$
	
	Using the standard basis $\{1, h\}$, we represent the hyperbolic number $z= x + hy$ as the vector $(x,y)\in \mathbb R^{1,1}.$ The real numbers $x$ and $y$ are called the real and unipotent parts of the hyperbolic number $z$ respectively. The hyperbolic conjugate $\bar{z}$ of $z$ is defined by $\bar{z}= x-hy$, and the hyperbolic modulus of $z$ is $$|z|_h=\sqrt{|\bar{z}z|}=\sqrt{|x^2-y^2|}.$$
	
	Let $z_1=x_1+hy_1, z_2=x_2+hy_2$, note that
	$$\bar{z}_1z_2=(x_1x_2-y_1y_2)+h(x_1y_2-x_2y_1),$$
	then the Minkowski inner product is
	$$\langle z_1,z_2\rangle_m=\Re {\bar{z_1}z_2}=x_1x_2-y_1y_2.$$
	In particular, 
	$$\bar{z}z=x^2-y^2=\langle z,z\rangle_m,$$
	so the hyperbolic modulus coincides with the Minkowski norm.
	
	The linear transformation $z\mapsto hz$ interchanges the components of $z=x+hy$, yielding $y+hx$. In the plane $\mathbb R^{1,1}$, this represents reflection symmetry about the line  $y=x.$ Moreover, $z$ and $hz$ are always orthogonal, just like $z$ and $iz$ in the complex plane. However, 
	$$\langle hz_1,hz_2\rangle_m=-\langle z_1,z_2\rangle_m.$$

	Though operationally similar to complex numbers, hyperbolic numbers fail to form a field due to zero divisors. The reason is that $(1-h)(1+h)=0$, so the null lines $y=\pm x$ (where$|z|_h =0$) consist entirely of non-invertible elements, reflecting their light-like character in $\mathbb R^{1,1}$.
	It is easy to verify that
	$$z^{-1}=\frac{\bar{z}}{\bar{z}z}=\frac{x-hy}{x^2-y^2}$$
	is the multiplicative inverse of $z$, whenever $|z|_h\ne 0$.

	The unit circle $\mathbb S^1\subset 
	\mathbb C$ is parameterized by $e^{i\theta}=\cos \theta + i\sin\theta$, where $\theta\in [0,2\pi)$ is the angle that the vector $z$ makes with the positive $x$-axis. Similarly, the right branch of the unit hyperbola $x^2-y^2=1$ is parameterized by $e^{h\theta}=\cosh \theta+h\sinh \theta $, where the parameter $\theta = \tanh^{-1}
	(y/x)\in (-\infty,+\infty)$ is called the \emph{hyperbolic angle}, or the rapidity of $z$.
	
	The set $\{e^{h\theta}|\theta\in \mathbb R\}$ forms a one-parameter group of hyperbolic rotations preserving $|z|_h$.  In relativistic mechanics, this group action models Lorentz boosts, with $\theta$ the rapidity of the boost.

	In addition to the standard basis $\{1, h\}$, which allows any hyperbolic number to be expressed as $z = x+hy$,
	the hyperbolic plane has another distinguished basis, associated with
	the two isotropic lines which separate the hyperbolic quadrants.
	The {\bf idempotent basis} of the hyperbolic numbers is the pair $\{h_+, h_-\}$,  where $h_+=\frac{1+h}{2}, h_-=\frac{1-h}{2}$.
	We say that $h_+$ and $h_-$ are idempotent, $h_+^2=h_+$ and $h_-^2=h_-$,
	and mutually annihilating, $h_+h_-=0$. In the basis $\{h_+, h_-\}$, the expression for $z=x+hy$ is $z=\hat{x}h_++\hat{y}h_-$, where
	$$\hat{x}=x+y, \hat{y}=x-y.$$
	
	For further background of the hyperbolic number plane, see \cite{hypernumber,MR2797703}
	
	\section{Curves in the Minkowski plane}\label{sec3}
	Let $\mathbf{\gamma}: I\to \mathbb{R}^{1,1}$ be a regular curve.
	At each point on the curve, the tangent vector $\mathbf{\gamma}_u$ is either space-like$(\langle \mathbf{\gamma}_u,\gamma_u\rangle_m >0)$, light-like$(\langle \gamma_u,\gamma_u\rangle_m =0)$, or time-like$(\langle \gamma_u,\gamma_u\rangle_m <0)$. The curvature $k(u)$ fails to exist when  $\gamma_u$ is light-like (due to degeneracy of the metric), so we will only look at curves where it is everywhere space-like or time-like to maintain $k(u)$ globally. 
	Since
	$$\langle h\gamma_u,h\gamma_u\rangle_m=-\langle \gamma_u,\gamma_u\rangle_m,$$
	then reflection
	across the $\hat{x}$-axis maps space-like curves to time-like and vice versa.
	
	The arc-length parameter $s$ for a curve $\mathbf{\gamma}: I\to \mathbb{R}^{1,1}$  is defined via the Minkowski metric: 
	$$ds=\sqrt{|\langle\gamma_u,\gamma_u\rangle_m|}du.$$ 
	The unit tangent vector is 
	$$\mathbf T=\gamma_s=\frac{1}{\sqrt{|\langle\gamma_u,\gamma_u\rangle_m|}}\gamma_u.$$
	The unit normal vector is obtained by reflecting $\mathbf T$ across the $\hat{x}$-axis (i.e. $\mathbf N=h\mathbf T$ ). Since $\langle\mathbf T,\mathbf T\rangle_m=\varepsilon$ (with $\varepsilon=\pm 1$ corresponds to space-like ($+1$) or time-like ($-1$) curves), the derivative $\mathbf T_s$ must be parallel to $\mathbf N$. The signed curvature $k$ is then defined by the structure equation $\mathbf T_s=k\mathbf N$.  Then we also have $\mathbf N_s=k\mathbf T$, so the Frenet formulas take the form 
	$$\mathbf T_s=k\mathbf N, \mathbf N_s=k\mathbf T.$$
	
	By definition, the mean curvature vector of $\gamma$ is given by $\mathbf H =\varepsilon\gamma_{ss}=\varepsilon\mathbf T_s=\varepsilon k\mathbf N $, so we always have $\langle \mathbf H, \mathbf N\rangle_m=-k$. 
	
	The map $z\mapsto hz$, that is the reflection across the $\hat{x}-$axis, reverses the sign of the Minkowski inner product and therefore transforms a time-like curve into a space-like one. Under it $\mathbf T\mapsto h\mathbf T$, $\mathbf N\mapsto\mathbf T$ and the curvature $k$ and the arc length $s$ are unchanged, so that the velocity $\frac{1}{k}\mathbf N$ is carried to $\frac{1}{k}h\mathbf N$; equivalently, both $\mathbf H$ and $\langle\mathbf H,\mathbf H\rangle_m$ change sign and the velocity $-\mathbf H/\langle\mathbf H,\mathbf H\rangle_m$ is preserved. Hence this reflection maps solutions of the inverse mean curvature flow to solutions, with the same time parameter, and without loss of generality we may restrict our analysis to space-like curves.
	Reversing the orientation if necessary, we will parametrize our space-like curves such that $\mathbf T$ lies in the right arm of the unit
	hyperbola $x^2-y^2=1$, then $\mathbf T=e^{h\theta}$, for some real number $\theta$ called the hyperbolic tangent angle. Note that $\mathbf T_s=\theta_she^{h\theta}=\theta_s\mathbf N$, so $\theta$ satisfies $\theta_s=k$.
	
	\subsection{Self-similar motions of Curves under the inverse mean curvature flow}
	
	Let $\mathbf{\gamma}: I\to \mathbb{R}^{1,1}$ be a curve parameterized by the arc length $s$. A self-similar motion of $\gamma$ is a map 
	$\hat{\mathbf{\gamma}}:I\times J\to \mathbb{R}^{1,1}$ of the form
	\begin{equation}
		\hat{\mathbf{\gamma}}(s,t)=g(t)e^{hf(t)}\mathbf{\gamma}(s)+\mathbf G(t).
	\end{equation}
	Here $J$ is an interval containing $0$ and $f,g:J\to \mathbb{R}$ and $\mathbf G: J\to \mathbb{R}^{1,1}$ are differentiable functions such that $f(0)=0, g(0)=1$ and $\mathbf G(0)=\mathbf 0$. Hence $\hat{\mathbf{\gamma}}(s,0)=\mathbf{\gamma}(s)$. The function $f$ determines the hyperbolic rotation, $g$ determines the dilation and $\mathbf G$ is the translation term.
	
	If $\hat{\mathbf{\gamma}}(s,t)$ is a solution of the inverse mean curvature flow, that is 
	$$\frac{\partial \hat{\mathbf{\gamma}}(s,t) }{\partial t}=\frac{1}{k(s,t)}\mathbf N(s,t),$$
	where $k(s,t)=\frac{k(s)}{g(t)}$ is the curvature and $\mathbf N(s,t)=e^{hf(t)}\mathbf N(s)$ is the normal vector of $\hat{\mathbf{\gamma}}(s,t)$  , then 
	\begin{equation}\label{cfeq}
		\langle\frac{\partial \hat{\mathbf{\gamma}}(s,t) }{\partial t}, \mathbf N(s,t) \rangle_m=-\frac{1}{k(s,t)},
	\end{equation}
	where we used $\langle\mathbf N,\mathbf N\rangle_m=-1$. As is customary for geometric flows, a solution is understood up to tangential reparametrization; accordingly \eqref{cfeq}, the equality of the normal velocities, is the equation to be solved.
	
	Since 
	\begin{equation*}
		\frac{\partial \hat{\mathbf{\gamma}}(s,t) }{\partial t}=(g'(t)e^{hf(t)}+hg(t)f'(t)e^{hf(t)})\mathbf{\gamma}(s)+\mathbf G'(t),
	\end{equation*}
	therefore,
	\begin{align*}
		-\frac{g(t)}{k(s)} &=-\frac{1}{k(s,t)}=\langle\frac{\partial \hat{\mathbf{\gamma}}(s,t) }{\partial t}, \mathbf N(s,t) \rangle_m\\ &=\langle(g'(t)e^{hf(t)}+hg(t)f'(t)e^{hf(t)})\mathbf{\gamma}(s)+\mathbf G'(t), e^{hf(t)}\mathbf N(s)\rangle_m\\
		&=g'(t)\langle\gamma(s), \mathbf N(s) \rangle_m+g(t)f'(t)\langle h\gamma(s), \mathbf N(s)\rangle_m+\langle e^{-hf(t)}\mathbf G'(t), \mathbf N(s)\rangle_m\\
		&=g'(t)\langle \gamma(s), \mathbf N(s)\rangle_m-g(t)f'(t)\langle \gamma(s), \mathbf T(s)\rangle_m+\langle e^{-hf(t)}\mathbf G'(t), \mathbf N(s)\rangle_m.
	\end{align*}
	We see $\gamma(s)$ has to satisfy at $t=0$
	\begin{equation}\label{meq}
		a\langle \gamma(s), \mathbf T(s)\rangle_m-b\langle \gamma(s), \mathbf N(s)\rangle_m-\langle C, \mathbf N(s)\rangle_m=\frac{1}{k(s)},
	\end{equation}
	where $a=f'(0), b=g'(0), C=\mathbf G'(0)$.
	
	Conversely, if the initial curve $\gamma$ satisfies equation (\ref{meq}), it will move in a self-similar manner under the inverse mean curvature flow. We will divide into three cases to prove this.
	
\paragraph{Translation} First assume $\gamma$ satisfies the pure translation equation
	
	\begin{equation}\label{tr}
		\langle C, \mathbf N(s)\rangle_m =-\frac{1}{k(s)}.   
	\end{equation}
	Take $\mathbf G(t)=Ct$, then $\hat{\mathbf{\gamma}}(s,t)=\gamma(s)+\mathbf G(t)$ satisfies equation (\ref{cfeq}). Thus, under the flow, the curve $\gamma$ translates with constant velocity vector $C.$
	
\paragraph{Dilation and hyperbolic rotation}
	If $\gamma$ satisfies the  equation
	\begin{equation}\label{dr}
		a\langle \gamma(s), \mathbf T(s)\rangle_m-b\langle \gamma(s), \mathbf N(s)\rangle_m=\frac{1}{k(s)}, 
	\end{equation}
	we can take $\hat{\mathbf{\gamma}}(s,t)=e^{(b+ha)t}\gamma(s)$. 
	Since
	\begin{align*}
		\langle\frac{\partial \hat{\mathbf{\gamma}}(s,t) }{\partial t}, \mathbf N(s,t) \rangle_m&=\langle (b+ha)e^{(b+ha)t}\gamma(s), e^{hat}\mathbf N(s)\rangle_m\\
		&=\langle (b+ha)e^{bt}\gamma(s), \mathbf N(s)\rangle_m\\
		&=e^{bt}(b\langle \gamma(s), \mathbf N(s)\rangle_m+\langle ha\gamma(s), \mathbf N(s)\rangle_m)\\
		&=e^{bt}(b\langle \gamma(s), \mathbf N(s)\rangle_m-a\langle \gamma(s), \mathbf T(s)\rangle_m)\\
		&=-\frac{e^{bt}}{k(s)}=-\frac{1}{k(s,t)},
	\end{align*}
	Therefore, the curve under the inverse mean curvature flow is $\hat{\mathbf{\gamma}}(s,t)=e^{(b+ha)t}\gamma(s)$.

	Under the base $\{h_+,h_-\}$, $e^{(b+ha)t}=e^{(b+a)t}h_++e^{(b-a)t}h_-$. If $a^2=b^2$, the factor is $e^{2at}h_++h_-(a=b)$ or $h_++e^{-2at}h_-(a=-b)$.  That is, the screw dilation reduces to a linear dilation of the one variable $h_+$ and $h_-$, leaving  the other one intact, when $a^2=b^2.$
	
	\paragraph{Dilation, hyperbolic rotation and translation} 
	
	Now, we assume $\gamma$ satisfies the full equation (\ref{meq}). We rewrite the equation as follows:
	\begin{equation*}
		-\langle (b+ha)\gamma(s)+C, \mathbf N(s)\rangle_m=\frac{1}{k(s)}.
	\end{equation*}
	If $a^2\ne b^2$, we can get rid of the translation term $C$ by translating $\gamma(s)$ by a fixed vector. If we put $\tilde{\gamma}(s)=\gamma(s)+\frac{C}{b+ha}$, then $\tilde{\gamma}(s)$ satisfies the equation of the form (\ref{dr}). Therefore, the solution $\hat{\gamma}$ moves by a combination of hyperbolic rotation and dilation with respect to the point $-\frac{C}{b+ha}.$
	
	If $a^2=b^2$, we let $f(t)=at, g(t)=e^{bt}$ as before.  The $\mathbf G(t)$ should satisfy $\mathbf G'(t)=Cg(t)e^{hf(t)}=Ce^{(b+ha)t}$. Set $C=c_1h_++c_2h_-$, then 
	
	\begin{equation*}
		\mathbf G'(t)=Ce^{(b+ha)t}=e^{(b+a)t}c_1h_++e^{(b-a)t}c_2h_-=\begin{cases}
			e^{2at}c_1h_++c_2h_-, \text{if } a=b\\
			c_1h_++e^{2bt}c_2h_-, \text{if } a=-b
		\end{cases} .
	\end{equation*}
	
	By integration, we have 
	\begin{equation*}
		\mathbf G(t)=\begin{cases}
			\frac{c_1}{2a}(e^{2at}-1)h_++c_2th_-, \text{if } a=b\\
			c_1th_++\frac{c_2}{2b}(e^{2bt}-1)h_-, \text{if } a=-b
		\end{cases} ,
	\end{equation*}
	where $a\ne0$ in the first case and $b\ne0$ in the second (if $a=b=0$ the motion is the pure translation already treated above). Letting $\hat{\gamma}(s,t)=e^{(b+ha)t}\gamma(s)+ \mathbf G(t)$, then $\hat{\gamma}(s,t)$ is a solution of the inverse mean curvature flow.
	
	\begin{thm}
		If $\hat{\gamma}:I\times J\to \mathbb R^{1,1}  $ is a self-similar solution of the inverse mean curvature flow, then $\gamma(s)=\hat{\gamma}(s,0)$ satisfies the equation (\ref{meq}). Conversely, if an initial curve $\gamma(s): I\to \mathbb R^{1,1}$ satisfies the equation (\ref{meq}), then it generates a self-similar solution of the inverse mean curvature flow in $\mathbb R^{1,1}$.
	\end{thm}

	Thus, if we want to find all self-similar solutions, we only need to find all initial curves $\gamma(s)$ that satisfy equation (\ref{meq}) with prescribed $a, b$ and $C$.
	
	Throughout the paper the solutions are described in terms of the hyperbolic tangent angle $\theta$, and the following terminology is used.
	
	\begin{defi}\label{def:complete}
		Let $\gamma$ be a solution of (\ref{meq}) written in terms of $\theta$ on a maximal interval $I\subseteq\mathbb{R}$, so that $\gamma_{\theta}=\frac{1}{k}e^{h\theta}$ and $ds=\bigl|\frac{1}{k}\bigr|d\theta$. A zero of $\frac{1}{k}$ is a point at which $\gamma_{\theta}$ vanishes and the curvature blows up; we call it a \emph{singular point} of $\gamma$, and we say that $\gamma$ is \emph{smooth} if $\frac{1}{k}$ has no zero on $I$, that is, if $\gamma$ is a regular space-like curve on all of $I$. We say that $\gamma$ is \emph{complete} if its arc length diverges at both ends of $I$.
	\end{defi}
	
	Since $\int\bigl|\frac{1}{k}\bigr|d\theta$ converges near an interior zero of $\frac{1}{k}$, each maximal regular arc of a singular solution has finite length on that side, and completeness in the above sense is a statement about the two ends of $I$ only. A solution may therefore be complete and singular at the same time, and the two properties are recorded separately below.
	
	\section{Translating solitons}\label{sec4}
	
	A curve $\gamma$ translates with constant velocity $C$ under the inverse mean curvature flow if and only if it satisfies the equation $\langle C, \mathbf N\rangle_m =-\frac{1}{k}$. In particular $C\ne0$, since $C=0$ would force $\frac{1}{k}=0$ and the curve would fail to be space-like.

	A homothety or an isometry of $\mathbb{R}^{1,1}$ carries such a curve to a translating soliton again, with the same time parameter. Indeed, under $\gamma\mapsto\lambda e^{h\theta_0}\gamma+v$ with $\lambda>0$ one has $\frac{1}{k}\mapsto\frac{\lambda}{k}$ while the unit normal is rotated, so that $C\mapsto\lambda e^{h\theta_0}C$; and the reflections across the $x$- and the $y$-axis replace $C$ by $\bar{C}$ and by $-\bar{C}$ respectively. These maps therefore do not leave the class of translating solitons, and they act on $\mathbb{R}^{1,1}\setminus\{0\}$ with exactly three orbits, namely the space-like, the time-like and the light-like vectors. Hence it suffices to consider the velocity vectors $1,\ h,\ 2h_+=(1,1)$. Let $\theta$ be the hyperbolic tangent angle, that is $\mathbf T=e^{h\theta}=\cosh{\theta}+h\sinh{\theta}$. Then $\frac{d\theta}{ds}=k,$ and the normal vector $\mathbf N=h\mathbf T=\sinh{\theta}+h\cosh{\theta}$.
	\paragraph{Time-like Case: $C=h=(0, 1)$.}
		In this case, $\frac{1}{k}=- \langle C, \mathbf N\rangle_m =\cosh{\theta}>0.$ Denote $\gamma(s)=x(s)+h y(s)$, then 
		$$\frac{dx}{d\theta}=\frac{dx}{ds}\times\frac{ds}{d\theta}=\cosh{\theta}\times \frac{1}{k}=\cosh^2{\theta}=\frac{1+\cosh{2\theta}}{2}, $$
		$$\frac{dy}{d\theta}=\frac{dy}{ds}\times\frac{ds}{d\theta}=\sinh{\theta}\times \frac{1}{k}=\sinh{\theta}\cosh{\theta}=\frac{\sinh{2\theta}}{2}. $$
		By integrating, we have
		\begin{equation}\label{try}
			\gamma=\frac{2\theta+\sinh{2\theta}}{4}+h\frac{\cosh{2\theta}}{4}.
		\end{equation}
		Since the arc length differential  is $ds=\frac{e^{\theta}+e^{-\theta}}{2}d\theta$, 
		 then this is a complete smooth curve.
			
		\paragraph{Space-like Case: $C=1=(1, 0).$ }	
		In this case,	$\frac{1}{k}=- \langle C, \mathbf N\rangle_m =-\sinh{\theta}.$ Then 
			$$\frac{dx}{d\theta}=\frac{dx}{ds}\times\frac{ds}{d\theta}=\cosh{\theta}\times \frac{1}{k}=-\sinh{\theta}\cosh{\theta}=-\frac{\sinh{2\theta}}{2}, $$
			$$\frac{dy}{d\theta}=\frac{dy}{ds}\times\frac{ds}{d\theta}=\sinh{\theta}\times \frac{1}{k}=-\sinh^2{\theta}=\frac{1-\cosh{2\theta}}{2}. $$
			By integrating, we have
			\begin{equation}
				\gamma=-\frac{\cosh{2\theta}}{4}+h\frac{2\theta-\sinh{2\theta}}{4}. 
			\end{equation}
			Since $k=-\frac{1}{\sinh{\theta}}$ and  the arc length differential  is  $ds=\frac{|e^{\theta}-e^{-\theta}|}{2}d\theta$,  
			then $\gamma$  is complete but has a singularity at $\theta =0$ .
		
		\paragraph{Light-like Case: $C=1+h =(1, 1).$}	
			In this case, $\frac{1}{k}=- \langle C, \mathbf N\rangle_m =\cosh{\theta}-\sinh{\theta}>0.$ Then 
			$$\frac{dx}{d\theta}=\cosh{\theta}\times \frac{1}{k}=\cosh{\theta}(\cosh{\theta}-\sinh{\theta})=\frac{1+\cosh{2\theta}-\sinh{2\theta}}{2}, $$
			$$\frac{dy}{d\theta}=\sinh{\theta}\times \frac{1}{k}=\sinh{\theta}(\cosh{\theta}-\sinh{\theta})=\frac{1-\cosh{2\theta}+\sinh{2\theta}}{2}. $$
			By integrating, we have
			\begin{equation}\label{trC(1,1)}
				\gamma=\frac{2\theta-e^{-2\theta}}{4}+h\frac{2\theta+e^{-2\theta}}{4},
			\end{equation}
		which also can be expressed as $x^2-y^2=-\frac{1}{2}\theta e^{-2\theta}.$
			
			Since the arc length differential  is $ds=e^{-\theta}d\theta$  and $\int_0^{+\infty}e^{-\theta}d\theta=1$ , then this curve is not complete.
		
	\begin{thm}\label{translating thm}
		Let $\gamma:I\to \mathbb R^{1,1}$ be a space-like curve in the Lorentz-Minkowski plane which evolves as a translating soliton of the inverse mean curvature flow with constant velocity vector $C\ne0$. Up to homotheties and isometries of $\mathbb{R}^{1,1}$, such curves are classified into three causal types according to $C$:
		\begin{itemize}
			\item Time-like Case ($C=(0,1)$): the solution is a complete and smooth curve, parameterized by the hyperbolic tangent angle $\theta$ as: $$\gamma=\frac{2\theta+\sinh{2\theta}}{4}+h\frac{\cosh{2\theta}}{4}.$$
			\item Space-like Case ($C=(1,0)$): the solution is a complete curve which is smooth except for a cusp singularity at $\theta=0$, given by:
				$$\gamma=-\frac{\cosh{2\theta}}{4}+h\frac{2\theta-\sinh{2\theta}}{4}.$$
			\item Light-like Case ($C=(1,1)$): the solution is smooth but not complete (it has finite arc length in one direction), expressed as:
					$$	\gamma=\frac{2\theta-e^{-2\theta}}{4}+h\frac{2\theta+e^{-2\theta}}{4}.$$
		\end{itemize}
	\end{thm}		
		
		Figure~\ref{fig:translators} illustrates the three cases.	
	\begin{figure}[htbp]
		\centering
		\begin{tikzpicture}
			\tikzset{asymp/.style={dashed,gray!60,line width=0.4pt}}
			
			\begin{axis}[width=4.5cm,height=5.5cm,at={(0,0)},anchor=south west,
				axis lines=middle,ticks=none,xmin=-3,xmax=3,ymin=-4,ymax=4,clip=true]
				\addplot[blue,thick,domain=-2:2,samples=100,variable=\t] 
				({\t/2 + sinh(2*\t)/4}, {cosh(2*\t)/4});
				\draw[asymp] (axis cs:0,0) -- (axis cs:2.8,2.8);
				\draw[asymp] (axis cs:0,0) -- (axis cs:-2.8,2.8);
				\node[font=\tiny,gray] at (axis cs:1.5,3.2) {$y=x$};
				\node[font=\tiny,gray] at (axis cs:-1.5,3.2) {$y=-x$};
			\end{axis}
			\node[anchor=north, font=\small] at ([yshift=-0.4cm]current axis.south) 
			{(i) Time-like Case};
			
			\begin{axis}[width=4.5cm,height=5.5cm,at={(5cm,0)},anchor=south west,
				axis lines=middle,ticks=none,xmin=-4,xmax=1,ymin=-4,ymax=4,clip=true]
				\addplot[blue,thick,domain=-2:2,samples=100,variable=\t] 
				({-cosh(2*\t)/4}, {\t/2 - sinh(2*\t)/4});
				\draw[asymp] (axis cs:0,0) -- (axis cs:-2.8,-2.8);
				\draw[asymp] (axis cs:0,0) -- (axis cs:-2.8,2.8);
				\node[font=\tiny,gray] at (axis cs:-3.2,1.5) {$y=-x$};
				\node[font=\tiny,gray] at (axis cs:-3.2,-1.5) {$y=x$};
				\addplot[only marks,mark=*,mark size=1.5pt,red] coordinates {(-0.25,0)};
				\node[font=\tiny,red,anchor=south west] at (axis cs:-0.2,0.2) {$\theta=0$};
			\end{axis}
			\node[anchor=north, font=\small] at ([yshift=-0.4cm]current axis.south) 
			{(ii) Space-like Case};
			
			\begin{axis}[width=4.5cm,height=5.5cm,at={(10cm,0)},anchor=south west,
				axis lines=middle,ticks=none,xmin=-3.5,xmax=1.5,ymin=-4,ymax=4,clip=true]
				\addplot[blue,thick,domain=-2:2,samples=100,variable=\t] 
				({\t/2 - exp(-2*\t)/4}, {\t/2 + exp(-2*\t)/4});
				\draw[asymp] (axis cs:0,0) -- (axis cs:1.2,1.2);
				\draw[asymp] (axis cs:0,0) -- (axis cs:-3,3);
				\node[font=\tiny,gray] at (axis cs:0.8,1.4) {$y=x$};
				\node[font=\tiny,gray] at (axis cs:-2.5,2.8) {$y=-x$};
			\end{axis}
			\node[anchor=north, font=\small] at ([yshift=-0.4cm]current axis.south) 
			{(iii) Light-like Case};
			
		\end{tikzpicture}
		\caption{Translating solitons in $\mathbb{R}^{1,1}$: (i) $C=h$, (ii) $C=1$, (iii) $C=1+h$. Cusps are marked in red; the dashed lines are the light-like directions $y=\pm x$.}
		\label{fig:translators}
	\end{figure}
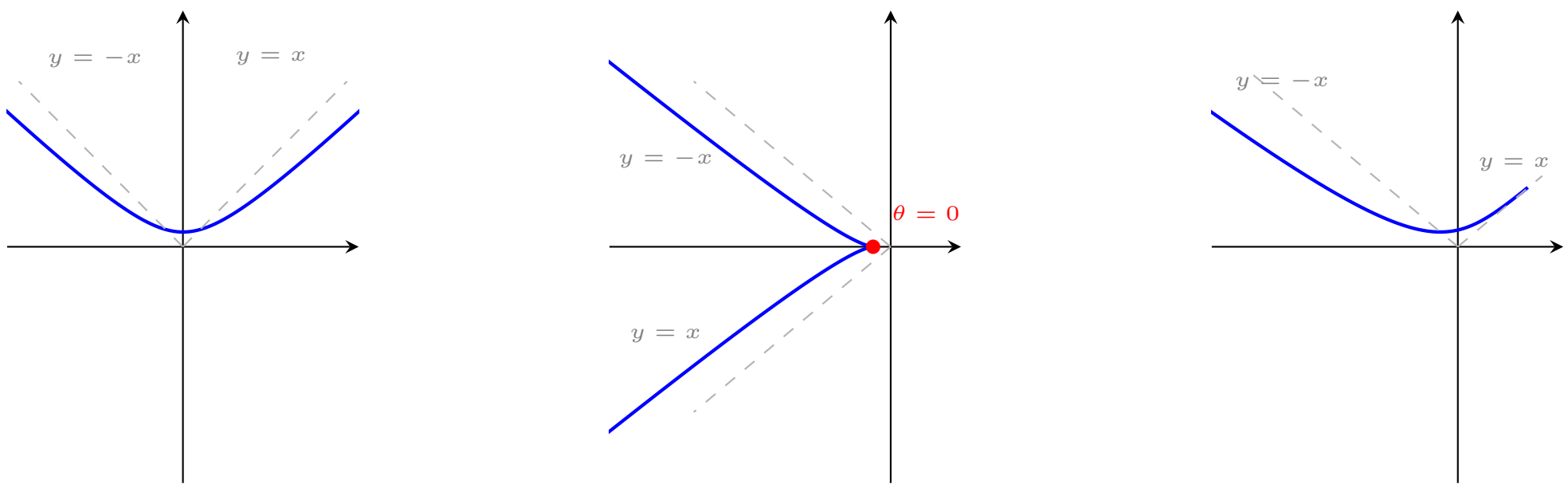
		\begin{rem}
			Formula \eqref{try} is the Wick rotation of the Euclidean translating soliton, the cycloid $(x,y)=(-\frac{\theta}{2}-\frac{\sin 2\theta}{4},\frac{\cos 2\theta}{4})$ of \cite{Drugan2016}: the substitutions $\theta\mapsto-i\theta$, $x\mapsto-ix$ together with $\cos iz=\cosh z$, $\sin iz=i\sinh z$ carry the cycloid to \eqref{try}.
		\end{rem}

		\section{Dilation and rotation}\label{sec5}
		
		From now on, we focus on curves satisfying the equation 
		
		\begin{equation}\label{eq}
			a\langle \gamma(s), \mathbf T(s)\rangle_m-b\langle \gamma(s), \mathbf N(s)\rangle_m=\frac{1}{k(s)}.
		\end{equation}
		Set $\tau(s)=\langle \gamma(s), \mathbf T(s)\rangle_m$ and $\nu(s)=\langle \gamma(s), \mathbf N(s)\rangle_m$, then $a\tau-b\nu=\frac{1}{k}$ and $\gamma=\tau\mathbf T-\nu \mathbf N$.
		By direct calculation, we have the following differential equations :
		\begin{equation}\label{tau}
			\begin{cases}
				\frac{d\tau}{ds}=1+k\nu=1+\frac{\nu}{a\tau-b\nu},\\
				\frac{d\nu}{ds}=k\tau=\frac{\tau}{a\tau-b\nu}.
			\end{cases}
		\end{equation}
		Note that whenever $a\tau-b\nu=0,$ the curvature $k$ is infinity and we have a singularity on the curve. We are interested in the smooth part.
		
		If we use the direction of the tangent vector $\theta$ as a parameter, 
		we have 
		$$\frac{d\theta}{ds}=k=\frac{1}{a\tau-b\nu}.$$
		Then the ODE system (\ref{tau}) about $\tau$ and $\nu$  can be rewritten as 
		\begin{equation}\label{systems of eq}
			\begin{cases}
				\frac{d\tau}{d\theta}=\frac{d\tau}{ds}\times \frac{ds}{d\theta}=a\tau+(1-b)\nu,\\
				\frac{d\nu}{d\theta}=\frac{d\nu}{ds}\times\frac{ds}{d\theta}=\tau.
			\end{cases} 
		\end{equation}
		This is a system of homogeneous linear differential equations with constant coefficients. We only need to find all the solutions to this system of equations (\ref{systems of eq}), and then, by considering the form $\gamma=\tau\mathbf T-\nu \mathbf N$, we can determine all the possible scenarios.
		
		First, we describe a special case, which $\tau(s)$ and $\nu(s)$ are linearly dependent along a curve.   In this case, the constant $a,b$ appearing in the relation (\ref{eq}) are not unique.
		
		\begin{pro}\label{hyper log spiral}
			$\tau(s)$ and $\nu(s)$ are linearly dependent along $\gamma$, say $\tau(s)=l\nu(s)$ with $l$ a constant, if and only if $l^2-al+b-1=0$. In this case $l\ne\pm1$ and, up to a homothety, $\gamma=(l-h)e^{(l+h)\theta}$.
		\end{pro}
		
		\begin{proof}
			If $\tau(s)$ and $\nu(s)$ are linearly dependent, we set $\tau(s)=l\nu(s)$ with $l$ a constant. then $a\tau(s)-b\nu(s)=(al-b)\nu(s)=\frac{1}{k(s)}$.
			By equation $(\ref{tau})$, we have 
			
			\begin{equation*}
				\begin{cases}
					\frac{d\tau}{ds}=l\frac{d\nu}{ds}=1+k\nu,\\
					\frac{d\nu}{ds}=k\tau=lk\nu,
				\end{cases}
			\end{equation*}
			then $\frac{1}{l}(1+k\nu)=lk\nu$.  Combining  $k=\frac{1}{(al-b)\nu}$ we have 
			
			\begin{equation}\label{te}
				l^2-al+b-1=0.
			\end{equation}
			Since $al-b=l^2-1$ by \eqref{te} and $(al-b)\nu=\frac{1}{k}\ne0$, we have $l\ne\pm1$.
			Let $T=e^{h\theta}$, then $N=hT=he^{h\theta}$ and $\gamma=\tau T-\nu N=(\tau-h\nu)e^{h\theta}=(l-h)\nu e^{h\theta}$.
			Since $\frac{d\nu}{d\theta}=\tau=l\nu$, then  $\nu=e^{l\theta}$ up to a constant and $\gamma=(l-h)e^{(l+h)\theta}.$
			
			Conversely, if $\gamma=(l-h)e^{(l+h)\theta}$ with $l\ne\pm 1$, then $\gamma_{\theta}=(l^2-1)e^{l\theta}e^{h\theta}$, a non-zero real multiple of $e^{h\theta}$.
			Hence, with the orientation for which $\frac{ds}{d\theta}=\frac{1}{k}=(l^2-1)e^{l\theta}$, we have $\mathbf T=\gamma_s=e^{h\theta}$, $\mathbf N=he^{h\theta}$ and $k=\frac{d\theta}{ds}=\frac{1}{(l^2-1)e^{l\theta}}$, so that $\theta$ is indeed the hyperbolic tangent angle of $\gamma$; note that $k<0$ when $l^2<1$, and that $ds=|l^2-1|e^{l\theta}d\theta$.
			$\tau=\langle \gamma, \mathbf T\rangle=le^{l\theta}$, $\nu=\langle \gamma, \mathbf N\rangle=e^{l\theta}$, that is $\tau=l\nu.$
			Since $l^2-al+b-1=0$, then $a\tau-b\nu=(al-b)\nu=(l^2-1)e^{l\theta}=\frac{1}{k}$.
		\end{proof}
		By \eqref{te}, $l=\pm1$ corresponds exactly to $b=\pm a$, and for these values $(l-h)e^{(l+h)\theta}$ is the constant $2h_-$, resp.\ $-2h_+$: the curve degenerates to a point, which is why $l=\pm1$ is excluded. For $l\ne\pm1$ one has $\langle\gamma,\gamma\rangle_m=(l^2-1)e^{2l\theta}$ and $ds=|l^2-1|e^{l\theta}d\theta$, so that in coordinates $(x,y)\in\mathbb R^{1,1}$ the curve satisfies
		$$x^2-y^2=(l^2-1)e^{2l\theta},\qquad \theta\in(-\infty,+\infty),$$
		with curvature radius $\bigl|\frac{1}{k}\bigr|=|l^2-1|e^{l\theta}$. Formally replacing the hyperbolic angle $\theta$ by an ordinary angle turns this relation into that of a logarithmic spiral $x^2+y^2=c\,e^{2l\theta}$ of $\mathbb R^2$; accordingly we call these curves \textbf{hyperbolic logarithmic spirals}. For $l=0$ the spiral is the standard hyperbola $y^2-x^2=1$, with constant curvature radius $1$, which is smooth and complete.
			Now, we begin to consider all solutions to the ODE system (\ref{systems of eq}). Let \begin{equation*} A=\left(
				\begin{matrix}
					a & 1-b\\
					1 & 0
				\end{matrix}\right),
			\end{equation*}
			the characteristic polynomial of $A$ is $|\lambda I-A|
				=\lambda^2-a\lambda+b-1,$  
				 which is same as equation (\ref{te}). 
				This is because if $\lambda$ is an eigenvalue of $A$, the corresponding eigenvector $(\tau, \nu)'$  satisfies $\tau=\lambda \nu$. Denote $\Delta=a^2-4(b-1)$. Next, we will analyze the different scenarios depending on the sign of $\Delta$.
				
			\subsection{Parabolic case: $\Delta=0$}
				
				In this case,  $A$ has a double eigenvalue $\lambda=\frac{a}{2}$, the corresponding eigenvector and generalized eigenvector are
				$v_1=\left(\begin{array}{c}
					\frac a 2\\
					1
				\end{array}\right)$ and $v_2=\left(\begin{array}{c}
					1\\
					0
				\end{array}\right)$ respectively.
				
				In this case, the general solution of the system of equations (\ref{systems of eq}) is 
				\begin{equation*}
					\left(\begin{matrix}
						\tau\\
						\nu
					\end{matrix}
					\right)=e^{\frac{a}{2}\theta}\left(C_1\left(\begin{matrix}
						\frac{a}{2}\\
						1
					\end{matrix}
					\right)+C_2\left(\begin{matrix}
						\frac{a}{2}\theta+1\\
						\theta
					\end{matrix}
					\right)\right)  , 
				\end{equation*}
				where $C_1, C_2$ are real constants determined by initial conditions that are not all zero.  From $\gamma=\tau\mathbf T-\nu \mathbf N$, the solution is given by 
				\begin{align*}
					\gamma&=\tau\mathbf T-\nu \mathbf N =\tau e^{h\theta}-h\nu e^{h\theta}=e^{h\theta}(\tau-h\nu)\\
					&=e^{(\frac{a}{2}+h)\theta}\left(C_1(\frac{a}{2}-h)+C_2(1+(\frac{a}{2}-h)\theta)\right).
				\end{align*}
				If $C_2=0,$ the curve can be expressed as $\gamma=(\frac{a}{2}-h) e^{(\frac{a}{2}+h)\theta}$,  which is the hyperbolic logarithmic spiral  described in proposition \ref{hyper log spiral} if $a\ne \pm 2$. If $a=\pm 2$, then $\frac{a}{2}-h$ is light-like and $\gamma$ degenerates to a point.
				
				If $C_2\ne 0$,  let $\tilde{\theta}=\theta+\frac{C_1}{C_2}$, then
				\begin{align*}
					\gamma&=e^{(\frac{a}{2}+h)\theta}\left(C_1(\frac{a}{2}-h)+C_2(1+(\frac{a}{2}-h)\theta)\right)\\
					&=e^{(\frac{a}{2}+h)(\tilde{\theta}-\frac{C_1}{C_2})}\left(C_1(\frac{a}{2}-h)+C_2(1+(\frac{a}{2}-h)(\tilde{\theta}-\frac{C_1}{C_2}))\right)\\
					&=C_2e^{-(\frac{a}{2}+h)\frac{C_1}{C_2}}e^{(\frac{a}{2}+h)\tilde{\theta}}(1+(\frac{a}{2}-h)\tilde{\theta}).
				\end{align*}
				That is, up to a homothety and a hyperbolic rotation, the curve can be expressed as 
				$$\gamma=e^{(\frac{a}{2}+h)\tilde{\theta}}(1+(\frac{a}{2}-h)\tilde{\theta}),$$
				where $\tilde{\theta}$  and $\theta$ differ by a constant. This  curve can be written as adding a hyperbolic logarithmic spiral $e^{(\frac{a}{2}+h)\tilde{\theta}}$ with $(\frac{a}{2}-h)\tilde{\theta}e^{(\frac{a}{2}+h)\tilde{\theta}} $  spiral. Since the curvature is $k=\bigl(a+(\tfrac{a^2}{4}-1)\tilde{\theta}\bigr)^{-1}e^{-\frac{a}{2}\tilde{\theta}}$, the curve is singular at $\tilde{\theta}=\frac{4a}{4-a^2}$, provided $a\ne\pm 2$. If $a=\pm 2$, then $b=2$ and $\frac{1}{k}=ae^{\pm\tilde{\theta}}$ never vanishes, so that the curve is smooth everywhere; in the light-cone coordinates $u=x+y$, $v=x-y$ it is the graph $u=e^{v-1}$ for $a=2$.
				By direct calculation, 
				$$ds=\sqrt{\langle\gamma_{\tilde{\theta}},\gamma_{\tilde{\theta}}\rangle_m}d\tilde{\theta}=|a+(\frac{a^2}{4}-1)\tilde{\theta}|e^{\frac{a}{2}\tilde{\theta}}d\tilde{\theta},$$
				then this curve is complete if and only if $a=0.$
				
				
				For $a=0$, $b=1$ the complete solution is $\gamma=(1-h\tilde{\theta})e^{h\tilde{\theta}}$, with a cusp at $\tilde{\theta}=0$ and the two asymptotes $y=\pm x$; see Figure~\ref{fig:dilrot}(i).
				\subsection{Hyperbolic case: $\Delta>0$}
				In this case, the matrix A has two different eigenvalues $\lambda_1, \lambda_2$. Let $K=\sqrt{\frac{a^2}{4}-b+1}$, then $\lambda_1=K+\frac{a}{2}$ and $\lambda_2=-K+\frac{a}{2}$ and the corresponding eigenvectors  are
				$v_1=\left(\begin{array}{c}
					K+\frac a 2\\
					1
				\end{array}\right)$ and 
				$v_2=\left(\begin{array}{c}
					-K+\frac a 2\\
					1
				\end{array}\right)$ respectively. 
				From the theory of systems of ordinary differential equations we know the general solution of (\ref{systems of eq}) is 
				\begin{equation*}
					\left(\begin{matrix}
						\tau\\
						\nu
					\end{matrix}
					\right)=C_1\left(\begin{matrix}
						\frac{a}{2}+K\\
						1
					\end{matrix}
					\right)e^{(\frac{a}{2}+K)\theta}+C_2\left(\begin{matrix}
						\frac{a}{2}-K\\
						1
					\end{matrix}
					\right)e^{(\frac{a}{2}-K)\theta}  , 
				\end{equation*}
				where $C_1, C_2$ are real constants that are not all zero.  
				From $\gamma=\tau\mathbf T-\nu \mathbf N$, the solution is given by 
				\begin{align*}
					\gamma&=\tau\mathbf T-\nu \mathbf N =\tau e^{h\theta}-h\nu e^{h\theta}=e^{h\theta}(\tau-h\nu)\\
					&=C_1(\frac{a}{2}+K-h)e^{(\frac{a}{2}+K+h)\theta}+C_2(\frac{a}{2}-K-h)e^{(\frac{a}{2}-K+h)\theta}\\
					&=C_1(\lambda_1-h)e^{(\lambda_1+h)\theta}+C_2(\lambda_2-h)e^{(\lambda_2+h)\theta}.
				\end{align*}
				If $C_2=0$,  we can write the solution as 
				$$\gamma=(\frac{a}{2}+K-h)e^{(\frac{a}{2}+K+h)\theta}=(\lambda_1-h)e^{(\lambda_1+h)\theta}.$$
				$\lambda_1^2=1$  if and only if $a=b<2$ or $-a=b>2$.
				
				If $C_1=0$, we can write the solution as 
				$$\gamma=(\frac{a}{2}-K-h)e^{(\frac{a}{2}-K+h)\theta}=(\lambda_2-h)e^{(\lambda_2+h)\theta}.$$
				$\lambda_2^2=1$  if and only if $a=b>2$ or $-a=b<2$.
				By Proposition \ref{hyper log spiral}, the above two curves are hyperbolic logarithmic spirals as long as $\lambda_i^2\ne 1$. On the two exceptional lines $b=\pm a$ one has $\lambda_i^2=1$ for exactly one index $i$ (both indices only when $a=b=0$), and then the corresponding curve degenerates to the point $2h_-$ or $-2h_+$ and has to be discarded. Moreover, since $\frac{1}{k}=(\lambda_i^2-1)e^{\lambda_i\theta}$ along such a spiral, its arc length converges at one of the two ends unless $\lambda_i=0$; thus none of these spirals is complete except for $\lambda_i=0$, in which case it is the standard hyperbola $y^2-x^2=1$.

				If $C_1\ne 0, C_2\ne 0$,  let  $\tilde{\theta}=\theta+\theta_0$,
				Since   
				\begin{align*}
					\gamma&=C_1(\lambda_1-h)e^{(\lambda_1+h)\theta}+C_2(\lambda_2-h)e^{(\lambda_2+h)\theta}\\
					&=C_1e^{-(\lambda_1+h)\theta_0}\left[(\lambda_1-h)e^{(\lambda_1+h)\tilde{\theta}}+\frac{C_2e^{-\lambda_2\theta_0}}{C_1e^{-\lambda_1\theta_0}}(\lambda_2-h)e^{(\lambda_2+h)\tilde{\theta}}\right]
				\end{align*}
				If $\frac{C_2}{C_1}>0$, we can choose $\theta_0$ such that $\frac{C_2e^{-\lambda_2\theta_0}}{C_1e^{-\lambda_1\theta_0}}=1$.   If $\frac{C_2}{C_1}<0$, we can choose $\theta_0$ such that $\frac{C_2e^{-\lambda_2\theta_0}}{C_1e^{-\lambda_1\theta_0}}=-1$.
				Then, up to a homothety and a hyperbolic rotation, these curves can be expressed as
			\begin{equation}\label{contraction}
				\gamma_1=(\lambda_1-h)e^{(\lambda_1+h)\tilde{\theta}}+(\lambda_2-h)e^{(\lambda_2+h)\tilde{\theta}},
			\end{equation}
				and 
					\begin{equation}\label{rota}
						\gamma_2=(\lambda_1-h)e^{(\lambda_1+h)\tilde{\theta}}-(\lambda_2-h)e^{(\lambda_2+h)\tilde{\theta}},
					 \end{equation}
				where $\tilde{\theta}$  and $\theta$ differ by a constant.
				
				Since
				$$\frac{d\gamma_1}{d\tilde{\theta}}=(\lambda_1^2-1)e^{(\lambda_1+h)\tilde{\theta}}+(\lambda_2^2-1)e^{(\lambda_2+h)\tilde{\theta}},
				$$
				then $\gamma_1$ has singularities if and only if $(\lambda_1^2-1)(\lambda_2^2-1)=b^2-a^2<0.$
				While $$\frac{d\gamma_2}{d\tilde{\theta}}=(\lambda_1^2-1)e^{(\lambda_1+h)\tilde{\theta}}-(\lambda_2^2-1)e^{(\lambda_2+h)\tilde{\theta}},
				$$
				then $\gamma_2$ has singularities if and only if $(\lambda_1^2-1)(\lambda_2^2-1)=b^2-a^2>0.$
				
				Since 
				$$ds_{1,2}=|(\lambda_1^2-1)e^{\lambda_1\tilde{\theta}}\pm(\lambda_2^2-1)e^{\lambda_2\tilde{\theta}}|d\tilde{\theta},$$
				the behaviour at the two ends is governed by those exponentials whose coefficients do not vanish. By \eqref{te}, $\lambda_i^2=1$ happens exactly when $b=\pm a$, and in that case the corresponding exponential disappears from $\frac{1}{k}$ altogether; we therefore distinguish three cases.
				\begin{itemize}
					\item If $b\ne\pm a$, both coefficients are non-zero, the end $\tilde{\theta}\to+\infty$ is dominated by $e^{\lambda_1\tilde{\theta}}$ and the end $\tilde{\theta}\to-\infty$ by $e^{\lambda_2\tilde{\theta}}$. Hence $\gamma_{1,2}$ is complete if and only if $\lambda_2\le 0\le\lambda_1$, that is, if and only if $\lambda_1\lambda_2=b-1\le 0$, or equivalently $b\le 1$.
					\item If exactly one of the two eigenvalues satisfies $\lambda_i^2=1$, only the other exponential survives in $\frac{1}{k}$, and both ends are governed by it; then $\gamma_{1,2}$ is complete if and only if the remaining eigenvalue vanishes. Since $\lambda_1+\lambda_2=a$ and $\lambda_1\lambda_2=b-1$, this occurs precisely for $a=\pm 1$ and $b=1$, in which case $\frac{1}{k}$ is constant.
					\item If $\lambda_1^2=\lambda_2^2=1$, then $\{\lambda_1,\lambda_2\}=\{1,-1\}$ and $a=b=0$, so that $\frac{1}{k}\equiv 0$ and there is no space-like solution at all. This is consistent with \eqref{eq}, which for $a=b=0$ reads $0=\frac{1}{k}$.
				\end{itemize}
				In particular the criterion $b\le 1$ is not sufficient on the two lines $b=\pm a$: for instance $a=b=-1$ gives $\lambda_1=1$, $\lambda_2=-2$ and $\frac{1}{k}=3e^{-2\tilde{\theta}}$, whose integral over $[0,+\infty)$ equals $\frac{3}{2}$, so that $\gamma_{1,2}$ is incomplete although $b\le 1$.
				
				\begin{rem}\label{rem:nulltrans}
					The lines $b=\pm a$ have an intrinsic meaning. For the light-like vectors $c=2h_{\mp}=(1,\mp1)$ one has $\langle c,\mathbf T\rangle_m=\pm\langle c,\mathbf N\rangle_m=e^{\pm\theta}$, so that $\gamma$ and $\gamma+c$ satisfy the same equation \eqref{eq} precisely when $a=\pm b$. On these lines \eqref{eq} thus acquires a translational symmetry along one light-like direction; correspondingly one of the spirals collapses to a point, and $\gamma_1$, $\gamma_2$ differ by such a translation.
				\end{rem}

				Three special cases deserve mention. For $a=b=1$, $\gamma_1=1-\sinh\tilde{\theta}-h(1+\cosh\tilde{\theta})$ is the standard hyperbola $(y+1)^2-(x-1)^2=1$. For a pure dilation $a=0$, the condition $\Delta>0$ forces $b<1$, so that $\gamma_{1,2}$ are complete and $\gamma_1$ is smooth; since $b=g'(0)$, all contracting solitons ($b<0$) belong to the hyperbolic case and are complete when $b\ne\pm a$. For a pure hyperbolic rotation $b=0$, again $\gamma_{1,2}$ are complete, now $\gamma_2$ is smooth, and the two spirals are incomplete since $\lambda_1\lambda_2=-1$. Figure~\ref{fig:dilrot}(ii),(iii) shows $\gamma_1,\gamma_2$ for $a=1$, $b=0$, where $\lambda_{1,2}=\frac{1\pm\sqrt{5}}{2}$.
			
			\subsection{Elliptic case: $\Delta<0$}

			In this case, the matrix A has two conjugate complex eigenvalues $\lambda, \bar{\lambda}$. Let $K=\sqrt{b-1-\frac{a^2}{4}}$, then $\lambda=iK+\frac{a}{2}$ and the corresponding eigenvectors  are
			$v_1=\left(\begin{array}{c}
				iK+\frac a 2\\
				1
			\end{array}\right)$ and 
			$v_2=\left(\begin{array}{c}
				-iK+\frac a 2\\
				1
			\end{array}\right)$ respectively.
			From the theory of systems of ordinary differential equations we know the general solution of (\ref{systems of eq}) is 
			\begin{align*}
				\left(\begin{matrix}
					\tau\\
					\nu
				\end{matrix}
				\right)&=C_1\left(\left(\begin{matrix}
					\frac{a}{2}\\
					1
				\end{matrix}
				\right)\cos{K\theta}-\left(\begin{matrix}
					K\\
					0
				\end{matrix}
				\right)\sin{K\theta}\right)e^{\frac{a}{2}\theta}
				\\&+C_2\left(\left(\begin{matrix}
					\frac{a}{2}\\
					1
				\end{matrix}
				\right)\sin{K\theta}+\left(\begin{matrix}
					K\\
					0
				\end{matrix}
				\right)\cos{K\theta}\right)e^{\frac{a}{2}\theta},
				\end{align*}
				where $C_1, C_2$ are constants that are not all zero.  
				From $\gamma=\tau\mathbf T-\nu \mathbf N$, the solution is given by 
				\begin{align*}
					\gamma&=\tau\mathbf T-\nu \mathbf N =\tau e^{h\theta}-h\nu e^{h\theta}=e^{h\theta}(\tau-h\nu)\\
					&=C_1e^{(\frac{a}{2}+h)\theta}(\frac{a}{2}\cos{K\theta}-K\sin{K\theta}-h\cos{K\theta})\\
					&+C_2e^{(\frac{a}{2}+h)\theta}(\frac{a}{2}\sin{K\theta}+K\cos{K\theta}-h\sin{K\theta})
				\end{align*}
				Up to a homothety and a hyperbolic rotation, we can write 
				$$\gamma=e^{(\frac{a}{2}+h)\tilde{\theta}}(\frac{a}{2}\sin{K\tilde{\theta}}+K\cos{K\tilde{\theta}}-h\sin{K\tilde{\theta}}),$$
				where $\tilde{\theta}$ and $\theta$ differ by a constant. Since 
				
				$$\gamma_{\tilde{\theta}}=e^{(\frac{a}{2}+h)\tilde{\theta}}((\frac{a^2}{2}-b)\sin{K\tilde{\theta}}+aK\cos{K\tilde{\theta}}),$$
				and 
				$$ds=|(\frac{a^2}{2}-b)\sin{K\tilde{\theta}}+aK\cos{K\tilde{\theta}}|e^{\frac{a}{2}\tilde{\theta}}d\tilde{\theta}, $$
				thus, this curve is complete if and only if $a=0,$  but it always has infinitely many singular points.
				
				For $a=0$ every solution of the elliptic case is complete, with infinitely many singular points; such a motion is an expansion, since $b>1$. The converse fails: an expansion with $0<b<1$ belongs to the hyperbolic case, where complete smooth solutions do exist. For $b=2$ the curve is $\gamma=e^{h\tilde{\theta}}(\cos\tilde{\theta}-h\sin\tilde{\theta})$, with cusps at $\tilde{\theta}\in\pi\mathbb Z$; see Figure~\ref{fig:dilrot}(iv).

				\begin{figure}[htbp]
					\centering
					\begin{tikzpicture}
						\pgfmathsetmacro{\phi}{(1+sqrt(5))/2}
						\pgfmathsetmacro{\psi}{(1-sqrt(5))/2}
						\tikzset{asymp/.style={dashed,gray!60,line width=0.4pt}}
						\begin{axis}[width=4.5cm,height=5.5cm,at={(0,0)},anchor=south west,
							axis lines=middle,ticks=none,xmin=-4,xmax=1.5,ymin=-4,ymax=4,clip=true]
							\addplot[blue,thick,domain=-2:2,samples=200,variable=\t]
							({cosh(\t) - \t*sinh(\t)}, {sinh(\t) - \t*cosh(\t)});
							\draw[asymp] (axis cs:0,0) -- (axis cs:-3.6,3.6);
							\draw[asymp] (axis cs:0,0) -- (axis cs:-3.6,-3.6);
							\node[font=\tiny,gray] at (axis cs:-2.4,3.4) {$y=-x$};
							\node[font=\tiny,gray] at (axis cs:-2.4,-3.4) {$y=x$};
							\addplot[only marks,mark=*,mark size=1.5pt,red] coordinates {(1,0)};
							\node[font=\tiny,red,anchor=south west] at (axis cs:0.6,0.3) {$\tilde{\theta}=0$};
						\end{axis}
						\node[anchor=north, font=\small] at ([yshift=-0.4cm]current axis.south)
						{(i) Parabolic: $a=0$, $b=1$};
						\begin{axis}[width=4.5cm,height=5.5cm,at={(5cm,0)},anchor=south west,
							axis lines=middle,ticks=none,xmin=-12,xmax=12,ymin=-12,ymax=12,clip=true]
							\addplot[blue,thick,domain=-2.5:1.2,samples=300,variable=\t]
							({exp(\phi*\t)*(\phi*cosh(\t)-sinh(\t)) + exp(\psi*\t)*(\psi*cosh(\t)-sinh(\t))},
							{exp(\phi*\t)*(\phi*sinh(\t)-cosh(\t)) + exp(\psi*\t)*(\psi*sinh(\t)-cosh(\t))});
							\draw[asymp] (axis cs:-11,11) -- (axis cs:11,-11);
							\draw[asymp] (axis cs:-11,-11) -- (axis cs:11,11);
							\node[font=\tiny,gray] at (axis cs:8.5,10.5) {$y=x$};
							\node[font=\tiny,gray] at (axis cs:-8.5,10.5) {$y=-x$};
							\addplot[only marks,mark=*,mark size=1.5pt,red] coordinates {(0.80,-1.97)};
						\end{axis}
						\node[anchor=north, font=\small] at ([yshift=-0.4cm]current axis.south)
						{(ii) Rotation: $\gamma_1$, $a=1$, $b=0$};
						\begin{axis}[width=4.5cm,height=5.5cm,at={(0,-6.3cm)},anchor=south west,
							axis lines=middle,ticks=none,xmin=-12,xmax=12,ymin=-12,ymax=12,clip=true]
							\addplot[blue,thick,domain=-2.5:1.2,samples=300,variable=\t]
							({exp(\phi*\t)*(\phi*cosh(\t)-sinh(\t)) - exp(\psi*\t)*(\psi*cosh(\t)-sinh(\t))},
							{exp(\phi*\t)*(\phi*sinh(\t)-cosh(\t)) - exp(\psi*\t)*(\psi*sinh(\t)-cosh(\t))});
							\draw[asymp] (axis cs:-11,11) -- (axis cs:11,-11);
							\draw[asymp] (axis cs:-11,-11) -- (axis cs:11,11);
							\node[font=\tiny,gray] at (axis cs:8.5,10.5) {$y=x$};
							\node[font=\tiny,gray] at (axis cs:-8.5,10.5) {$y=-x$};
						\end{axis}
						\node[anchor=north, font=\small] at ([yshift=-0.4cm]current axis.south)
						{(iii) Rotation: $\gamma_2$, $a=1$, $b=0$};
						\begin{axis}[width=4.5cm,height=5.5cm,at={(5cm,-6.3cm)},anchor=south west,
							axis lines=middle,ticks=none,xmin=-13,xmax=4,ymin=-13,ymax=13,clip=true]
							\addplot[blue,thick,domain=-4:4,samples=500,variable=\t]
							({cosh(\t)*cos(\t r) - sinh(\t)*sin(\t r)},
							{sinh(\t)*cos(\t r) - cosh(\t)*sin(\t r)});
							\draw[asymp] (axis cs:0,0) -- (axis cs:-12.5,12.5);
							\draw[asymp] (axis cs:0,0) -- (axis cs:-12.5,-12.5);
							\node[font=\tiny,gray] at (axis cs:-8.5,12) {$y=-x$};
							\node[font=\tiny,gray] at (axis cs:-8.5,-12) {$y=x$};
							\addplot[only marks,mark=*,mark size=1.5pt,red] coordinates {(1,0) (-11.59,11.55) (-11.59,-11.55)};
						\end{axis}
						\node[anchor=north, font=\small] at ([yshift=-0.4cm]current axis.south)
						{(iv) Elliptic: $a=0$, $b=2$};
					\end{tikzpicture}
					\caption{The three cases of Theorem~\ref{main thm1}: (i) parabolic, $a=0$, $b=1$; (ii), (iii) hyperbolic, $a=1$, $b=0$; (iv) elliptic, $a=0$, $b=2$, $\tilde{\theta}\in[-4,4]$. Cusps are marked in red; the dashed lines are the light-like directions $y=\pm x$.}
					\label{fig:dilrot}
				\end{figure}
\subsection{Summary}			
\begin{thm}	\label{main thm1}
	Under the invariance of scaling and rotation, the solutions can be classified into the following three cases:
	\begin{itemize}
		\item \textbf{Elliptic case: $\Delta<0$.}
		
		 For given  values of $a$ and $b$, the solution is unique in this case, and can be expressed as 
		$$\gamma=e^{(\frac{a}{2}+h)\tilde{\theta}}(\frac{a}{2}\sin{K\tilde{\theta}}+K\cos{K\tilde{\theta}}-h\sin{K\tilde{\theta}}),$$
		where $K=\sqrt{b-1-\frac{a^2}{4}}$ and $\tilde{\theta}$ and tangent angle $\theta$ differ by a constant. In this case, $\gamma$ is complete if and only if $a=0$ and has infinitely many singular points for all possible $a$. 
		
		\item \textbf{Parabolic case: $\Delta=0$.}	
		
		For given  $a$ and $b$, there are two distinct types of solutions in this case. One is the hyperbolic logarithmic spiral $\gamma=(\frac{a}{2}-h) e^{(\frac{a}{2}+h)\theta},$  and the other is 
		$$\gamma=e^{(\frac{a}{2}+h)\tilde{\theta}}(1+(\frac{a}{2}-h)\tilde{\theta}),$$
		where $\tilde{\theta}$  and  $\theta$ differ by a constant. The latter $\gamma$ is complete if and only if $a=0$; it has a cusp singularity at $\tilde{\theta}=\frac{4a}{4-a^2}$ when $a\ne\pm2$, and is smooth when $a=\pm 2$. For $a=\pm2$ the former $\gamma$ degenerates to a point.
		
	\item \textbf{Hyperbolic case: $\Delta>0$.}	
		
	Assume $(a,b)\ne(0,0)$. For given $a$ and $b$ with $b\ne\pm a$, there are four distinct types of solutions in this case. Two of them are hyperbolic logarithmic spirals, which we denote by 
	$(\lambda_1-h)e^{(\lambda_1+h)\theta}$ and 
	$(\lambda_2-h)e^{(\lambda_2+h)\theta}$, respectively. 
	The other two curves, denoted as	 	
	 \begin{equation*}
	 	\gamma_{1,2}=(\lambda_1-h)e^{(\lambda_1+h)\tilde{\theta}}\pm(\lambda_2-h)e^{(\lambda_2+h)\tilde{\theta}},
 	\end{equation*}
	where $K=\sqrt{\frac{a^2}{4}-b+1}$, and  $\lambda_1=K+\frac{a}{2}$, $\lambda_2=-K+\frac{a}{2}$. In this case, $\gamma_1$ is smooth if and only if $b^2\ge a^2$, while $\gamma_2$ is smooth if and only if $b^2\le a^2$; in particular both are smooth when $b^2=a^2$. Moreover $\gamma_1,\gamma_2$ are complete if and only if $b\le 1$.
	
	On the two exceptional lines $b=\pm a$ one of the eigenvalues satisfies $\lambda_i^2=1$; the corresponding spiral degenerates to a point, $\gamma_1$ and $\gamma_2$ differ by a translation, and they are complete if and only if $a=\pm 1$ and $b=1$. For $a=b=0$ there is no space-like solution.
	\end{itemize}	
		\end{thm}	
	
	\begin{cor} In particular, when we restrict our consideration to invariance under hyperbolic rotations or scaling alone, we obtain the following conclusions:
	\begin{itemize} 
	\item 	If $b=0$ and $a\ne 0$, so that the motion is a pure hyperbolic rotation, then $b<1$ and both solutions $\gamma_1,\gamma_2$ of \eqref{contraction} and \eqref{rota} are complete. Moreover $\gamma_2$ is smooth, so that for every non-zero rotation parameter $a$ there exists a solution that is both complete and smooth. The two hyperbolic logarithmic spirals belonging to the same $(a,b)$ are, on the contrary, never complete, since $\lambda_1\lambda_2=b-1=-1$ and hence neither exponent vanishes.
	\item  If $a=0$ and $b\ne 0$, so that the motion is a pure dilation, then the hyperbolic case $\Delta>0$ forces $b<1$, and again $\gamma_1,\gamma_2$ are complete while the two spirals are not. In addition $\gamma_1$ is smooth, so that for every dilation parameter $b<1$, $b\ne0$, there exists a solution that is both complete and smooth, namely $\gamma_1$ defined by equation $(\ref{contraction}).$ In particular this applies to every contracting pure dilation, i.e.\ to every $b<0$.
	\item If $a=1,b=1$, $\gamma_1$ degenerates to a standard hyperbola  $(y+1)^2-(x-1)^2=1.$ 
	\end{itemize}
	\end{cor}

\section{Dilation, rotation and translation}\label{sec6}
				
				Finally, we consider the whole equation (\ref{meq}) with $a^2=b^2$:
				\begin{equation*}
					a\langle \gamma(s), \mathbf T(s)\rangle_m-b\langle \gamma(s), \mathbf N(s)\rangle_m-\langle C, \mathbf N(s)\rangle_m=\frac{1}{k(s)}.
				\end{equation*}
				Since $a^{2}=b^{2}$, either $a=b$ or $a=-b$. The reflection in the $y$-axis, $\tilde{\gamma}(s)=-\overline{\gamma(-s)}$, is an isometry of $\mathbb{R}^{1,1}$ which preserves space-like curves and carries the data $(a,b,C)$ to $(-a,b,-\bar{C})$; in particular it leaves $b$ unchanged. If $C$ is space-like or time-like, then $-\bar{C}$ is again a multiple of the same normalized vector, so this reflection carries $a=-b$ to $a=b$ and we may assume $a=b$. If $C$ is light-like, the same reflection interchanges the two null directions $1+h$ and $1-h$; we therefore keep $a=b$ and treat the two light-like translations $C=(1,1)$ and $C=(1,-1)$ separately.

				We stress that the \emph{sign} of $a=b$ cannot be normalized away. For the inverse mean curvature flow the rescaling $\gamma\mapsto\lambda\gamma$ $(\lambda>0)$ gives $\frac{1}{k}\mapsto\frac{\lambda}{k}$, so it leaves $(a,b)$ invariant and only replaces $C$ by $\lambda C$; thus $a$ and $b$ are genuine moduli. This is in contrast with the mean curvature flow, for which the same rescaling acts by $(a,b)\mapsto(a/\lambda^{2},b/\lambda^{2})$ and permits the normalization $|b|=1$ \cite{Halldorsson2015}. Accordingly, in this section $a=b$ denotes an arbitrary real number, and we exclude $a=0$, which is the pure translation case of Section \ref{sec4}.
				Then take $\frac{1}{k}=a(\tau-\nu)-\langle C, \mathbf N\rangle_m$ into (\ref{systems of eq}), we get
				
				\begin{equation}
					\begin{cases}
						\frac{d\tau}{d\theta}=\frac{d\tau}{ds}\times \frac{ds}{d\theta}=\frac{(1+k\nu)}{k}=\frac{1}{k}+\nu=a\tau+(1-a)\nu-\langle C, \mathbf N\rangle_m,\\
						\frac{d\nu}{d\theta}=\frac{d\nu}{ds}\times\frac{ds}{d\theta}=\tau.
					\end{cases} 
				\end{equation}
				Then we have
					\begin{equation*} \left(
					\begin{matrix}
						\frac{d\tau}{d\theta}\\
						\frac{d\nu}{d\theta}
					\end{matrix}\right)=A\left(\begin{matrix}
						\tau\\
						\nu
					\end{matrix}\right)+f(\theta),
				\end{equation*}
				
			where
				\begin{equation*} A=\left(
					\begin{matrix}
						a & 1-a\\
						1 & 0
					\end{matrix}\right),
				\end{equation*}
			and  
				\begin{equation*} f(\theta)=\left(
					\begin{matrix}
						-\langle C, \mathbf N\rangle_m\\
						0
					\end{matrix}\right).
				\end{equation*}

				The characteristic polynomial  is  
				$$|\lambda I-A|=(\lambda-\frac{a}{2})^2-\frac{(a-2)^2}{4}=(\lambda-1)\bigl(\lambda-(a-1)\bigr),$$
				so that $a=2$ is precisely the resonant case in which the two eigenvalues coincide.
				
				It is convenient to record the equation in the light-cone coordinates
				$$\hat{x}=x+y,\qquad \hat{y}=x-y,\qquad \gamma=\hat{x}h_++\hat{y}h_-,$$
				and to write $C_{\pm}=C_x\pm C_y$ for the corresponding components of the translation vector. Since $\gamma_{\theta}=\frac{1}{k}e^{h\theta}$ we have
				\begin{equation}\label{eq:nullder}
					\frac{d\hat{x}}{d\theta}=\frac{1}{k}e^{\theta},\qquad \frac{d\hat{y}}{d\theta}=\frac{1}{k}e^{-\theta},
				\end{equation}
				while a direct computation gives $\tau-\nu=\hat{x}e^{-\theta}$ and $\langle C,\mathbf N\rangle_m=\frac{1}{2}\bigl(C_-e^{\theta}-C_+e^{-\theta}\bigr)$. Substituting these into $\frac{1}{k}=a(\tau-\nu)-\langle C,\mathbf N\rangle_m$ and using the first equation of \eqref{eq:nullder}, the whole system collapses to a single first order linear equation for $\hat{x}$:
				\begin{equation}\label{eq:xhat}
					\frac{d\hat{x}}{d\theta}=a\hat{x}+\frac{C_+}{2}-\frac{C_-}{2}e^{2\theta}.
				\end{equation}
				Once $\hat{x}$ is known, $\frac{1}{k}=\hat{x}'e^{-\theta}$ and $\hat{y}$ is obtained from \eqref{eq:nullder} by one further integration. We shall use \eqref{eq:xhat} together with the following elementary lemma.
				
				\begin{lem}\label{lem:normalisation}
					Let $a=b\ne0$ and let $\gamma$ satisfy \eqref{meq}. Then:
					\begin{enumerate}
						\item[(i)] a translation $\gamma\mapsto\gamma+v$ adds constants to $\hat{x}$ and $\hat{y}$, leaves $C_-$ and $\frac1k$ unchanged, and replaces $C_+$ by $C_+-2a(v_1+v_2)$;
						\item[(ii)] a homothety and a hyperbolic rotation $\gamma\mapsto\lambda e^{h\phi}\gamma$ with $\lambda>0$ replace $\theta$ by $\theta+\phi$ and multiply $C_-$ by $\lambda e^{-\phi}$, $C_+$ by $\lambda e^{\phi}$ and $\frac{1}{k}$ by $\lambda$;
						\item[(iii)] the map $\gamma(s)\mapsto-\gamma(-s)$ replaces $(C_+,C_-,\frac1k)$ by $(-C_+,-C_-,-\frac1k)$;
						\item[(iv)] two solutions with the same $a=b$ are congruent under a composition of the maps in (i)--(iii) if and only if their functions $\frac{1}{k}(\theta)$ agree after a shift $\theta\mapsto\theta+\phi$ and multiplication by a non-zero constant.
					\end{enumerate}
					In particular the translation vector is only defined modulo the light-like direction $(1,1)$: by (i) one may always normalize $C_+=0$, so that $C=\frac{C_-}{2}(1,-1)$, and the classification is governed by whether $C_-$ vanishes, and not by the causal character of $C$ itself.
				\end{lem}

				\begin{proof}
					(i) A translation does not change $\gamma_{\theta}$, hence neither $\mathbf T$, $\mathbf N$, $\theta$ nor $\frac{1}{k}$, while $\hat{x}\mapsto\hat{x}+w$ with $w=v_1+v_2$. Substituting into \eqref{eq:xhat} gives $\hat{x}'=a(\hat{x}+w)+\frac{C_+-2aw}{2}-\frac{C_-}{2}e^{2\theta}$, that is, $C_+\mapsto C_+-2aw$ while $C_-$ is unchanged.

					(ii) Since the boost $e^{h\phi}$ is an isometry, $\tilde{\gamma}=\lambda e^{h\phi}\gamma$ has unit tangent $e^{h\phi}\mathbf T$, hence tangent angle $\theta+\phi$, and $\tilde{\tau}=\lambda\tau$, $\tilde{\nu}=\lambda\nu$, $\frac{1}{\tilde{k}}=\frac{\lambda}{k}$. Inserting this into \eqref{meq} yields $\langle\tilde{C},e^{h\phi}\mathbf N\rangle_m=\lambda\langle C,\mathbf N\rangle_m$ for all $\theta$; as the vectors $\mathbf N(\theta)$ span $\mathbb R^{1,1}$, this means $\tilde{C}=\lambda e^{h\phi}C$. In the idempotent basis $e^{h\phi}=e^{\phi}h_++e^{-\phi}h_-$ and $C=C_+h_++C_-h_-$, so that $\tilde{C}_{\pm}=\lambda e^{\pm\phi}C_{\pm}$.

					(iii) For $\tilde{\gamma}(s)=-\gamma(-s)$ one has $\tilde{\mathbf T}(s)=\mathbf T(-s)$, so the orientation convention of Section \ref{sec3} is preserved and $\tilde{\theta}(s)=\theta(-s)$, while $\tilde{k}(s)=-k(-s)$, $\tilde{\tau}(s)=-\tau(-s)$ and $\tilde{\nu}(s)=-\nu(-s)$. Both sides of \eqref{meq} change sign, whence $\tilde{C}=-C$, and as functions of $\theta$ we get $\frac{1}{\tilde{k}}=-\frac{1}{k}$.

					(iv) Since $\gamma_{\theta}=\frac{1}{k}e^{h\theta}$, the function $\frac{1}{k}(\theta)$ determines $\gamma$ up to a translation. If $\frac{1}{\tilde{k}}(\theta)=c\,\frac{1}{k}(\theta+\phi)$ with $c\ne0$, apply (ii) to $\gamma$ with $\lambda=|c|$ and rotation angle $-\phi$, and (iii) if $c<0$; the resulting curve has the same function $\frac{1}{k}$ as $\tilde{\gamma}$ and therefore differs from it by a translation. The converse is immediate from (i)--(iii). We note that the maps (i)--(iii) exhaust the congruences within the class $a=b\ne0$: the homotheties and isometries of $\mathbb R^{1,1}$ are generated by translations, dilations $\lambda>0$, boosts, $-\mathrm{id}$ and the reflection $\gamma\mapsto\bar{\gamma}$, and the latter carries the data $(a,b,C)$ to $(-a,b,\bar{C})$, hence out of the class.
				\end{proof}
				
				 Next, we will classify and discuss the solutions based on the eigenvalues. 
				\subsection{The case $a=2$}
				
				In this case the two eigenvalues of $A$ coincide and the inhomogeneous term $e^{2\theta}$ of \eqref{eq:xhat} is resonant. The general solution of \eqref{eq:xhat} is
				$$\hat{x}=\Bigl(K-\frac{C_-}{2}\theta\Bigr)e^{2\theta}-\frac{C_+}{4},\qquad K\in\mathbb{R},$$
				whence, by \eqref{eq:nullder},
				$$\frac{1}{k}=\Bigl(2K-\frac{C_-}{2}-C_-\theta\Bigr)e^{\theta},\qquad
				\hat{y}=\Bigl(2K-\frac{C_-}{2}\Bigr)\theta-\frac{C_-}{2}\theta^{2}+D,\qquad D\in\mathbb{R}.$$
				By Lemma \ref{lem:normalisation}(i) the constants $C_+$ and $D$ are pure translations and may be normalized to zero, and by Lemma \ref{lem:normalisation}(iv) the congruence class of the solution can be read off from $\frac{1}{k}$ up to a shift of $\theta$ and a non-zero factor. Two cases occur.
				
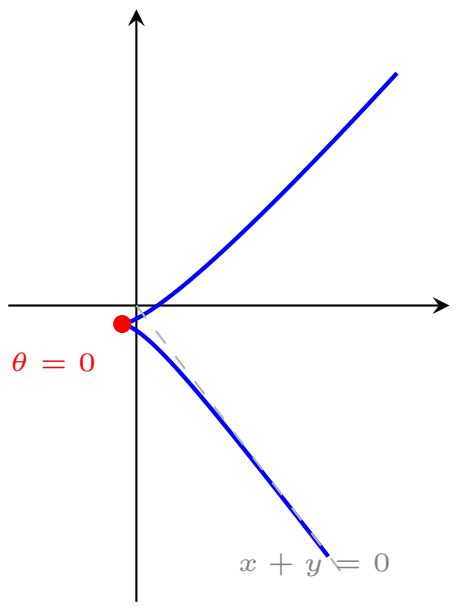
\begin{figure}[htbp]
					\centering
					\begin{tikzpicture}[
						declare function={
							gx(\t) = (2*\t-1)*exp(2*\t) + 2*\t^2;
							gy(\t) = (2*\t-1)*exp(2*\t) - 2*\t^2;
						}
						]
						\tikzset{asymp/.style={dashed,gray!60,line width=0.4pt}}
						\begin{axis}[width=4.5cm,height=5.5cm,at={(0,0)},anchor=south west,
							axis lines=middle,ticks=none,
							xmin=-9,xmax=22,ymin=-16,ymax=16,clip=true]
							\addplot[blue,thick,domain=-2.6:1.2,samples=300,variable=\t]
							({gx(\t)}, {gy(\t)});
							\draw[asymp] (axis cs:0,0) -- (axis cs:15,-15);
							\node[font=\tiny,gray] at (axis cs:12.5,-14) {$x+y=0$};
							\addplot[only marks,mark=*,mark size=1.5pt,red] coordinates {(-1,-1)};
							\node[font=\tiny,red,anchor=north east] at (axis cs:-1.5,-1.5) {$\theta=0$};
						\end{axis}
					\end{tikzpicture}
					\caption{The mixed soliton \eqref{eq:a2normal} ($a=b=2$, $C_-\ne0$), with its cusp at $\theta=0$ and the asymptote $x+y=0$.}
					\label{fig:mixeda2}
				\end{figure}

				\paragraph{The generic case $C_-\ne0$}
				Here $\frac{1}{k}=-C_-(\theta-\theta_1)e^{\theta}$ with $\theta_1=\frac{2K}{C_-}-\frac{1}{2}$, so that $K$ is absorbed by the shift $\theta\mapsto\theta+\theta_1$. After that shift and a homothety we may assume $\frac{1}{k}=8\theta e^{\theta}$, and then $\hat{x}=2(2\theta-1)e^{2\theta}$ and $\hat{y}=4\theta^{2}$; that is,
				\begin{equation}\label{eq:a2normal}
					\gamma=(2\theta-1)e^{2\theta}+2\theta^{2}
					+h\bigl((2\theta-1)e^{2\theta}-2\theta^{2}\bigr).
				\end{equation}
				The arc length differential formula is $ds=8|\theta|e^{\theta}d\theta$. Consequently $\gamma$ has a cusp singularity at $\theta=0$, at the point $(-1,-1)$; as $\theta\to-\infty$ one has $\hat{x}\to0$ and $\hat{y}\to+\infty$, so that $\gamma$ is asymptotic to the light-like line $x+y=0$; and its arc length is finite as $\theta\to-\infty$ and infinite as $\theta\to+\infty$, so that $\gamma$ is not complete. See Figure~\ref{fig:mixeda2}.
				
				In particular, the three normalized velocities $C=(1,0)$, $C=(0,1)$ and $C=(1,-1)$ all satisfy $C_-\ne0$ and therefore produce one and the same curve up to a homothety and an isometry: they give
				$$\frac{1}{k}=\Bigl(2K-\frac{1}{2}-\theta\Bigr)e^{\theta},\qquad
				\Bigl(2K+\frac{1}{2}+\theta\Bigr)e^{\theta},\qquad
				\bigl(2K-1-2\theta\bigr)e^{\theta}$$
				respectively, and each of these is a non-zero multiple of $\theta e^{\theta}$ after a shift of $\theta$. This is an instance of Lemma \ref{lem:normalisation}: modulo the light-like direction $(1,1)$, each of the three vectors is a multiple of $(1,-1)$.
				
				\paragraph{The degenerate case $C_-=0$}
				Here $C$ is parallel to the light-like direction $(1,1)$ and, by Lemma \ref{lem:normalisation}(i), a translation removes the translation term altogether: the equation reduces to \eqref{dr} with $a=b=2$, that is, to the parabolic case of Section \ref{sec5}. Explicitly $\frac{1}{k}=2Ke^{\theta}$, which never vanishes when $K\ne0$, and up to a homothety and a translation
				\begin{equation}\label{log sprial}
					\gamma=e^{2\theta}+2\theta+h(e^{2\theta}-2\theta),
				\end{equation}
				which is the light-like translating soliton $(\ref{trC(1,1)})$, a smooth but incomplete curve. If $K=0$ then $\frac{1}{k}\equiv0$ and there is no space-like solution.

				\subsection{The case $a\ne2$}
				
				Assume now $a=b\notin\{0,2\}$. The general solution of \eqref{eq:xhat} is
				$$\hat{x}=Ke^{a\theta}+\frac{C_-}{2(a-2)}e^{2\theta}-\frac{C_+}{2a},\qquad K\in\mathbb{R},$$
				and \eqref{eq:nullder} gives
				\begin{equation}\label{eq:kab}
					\frac{1}{k}=\alpha\,e^{(a-1)\theta}+\beta\,e^{\theta},\qquad \alpha=aK,\quad \beta=\frac{C_-}{a-2},
				\end{equation}
				$$\hat{y}=\frac{\alpha}{a-2}e^{(a-2)\theta}+\beta\theta+D,\qquad D\in\mathbb{R}.$$
				As before, $C_+$ and $D$ are removed by a translation. By Lemma \ref{lem:normalisation}(ii) a homothety and a hyperbolic rotation multiply $\alpha$ and $\beta$ by $\lambda e^{(a-1)\phi}$ and $\lambda e^{\phi}$ respectively; since $a\ne2$ these two factors can be prescribed independently, so that whenever $\alpha\beta\ne0$ we may normalize $|\alpha|=|\beta|=1$. Finally Lemma \ref{lem:normalisation}(iii) changes the signs of $\alpha$ and $\beta$ simultaneously. Hence, in contrast with the resonant case $a=2$, the constant $K$ can no longer be absorbed: the sign of $\alpha\beta$ is a congruence invariant, and by Lemma \ref{lem:normalisation}(iv) there are exactly four congruence classes of solutions for every value of $a$, namely
				$$\text{(I) } \alpha\beta>0,\qquad \text{(II) } \alpha\beta<0,\qquad \text{(III) } \alpha=0\ne\beta,\qquad \text{(IV) } \beta=0\ne\alpha,$$
				while $\alpha=\beta=0$ gives $\frac{1}{k}\equiv0$ and no space-like curve. The picture is entirely parallel to the hyperbolic case of Section \ref{sec5}, where the four types were the two hyperbolic logarithmic spirals and the two curves $\gamma_1,\gamma_2$ of \eqref{contraction} and \eqref{rota}.
				
\begin{figure}[htbp]
					\centering
					\begin{tikzpicture}[
						declare function={
							px(\a,\s,\t) = 0.5*( exp(\a*\t)/\a + \s*exp(2*\t)/2 + exp((\a-2)*\t)/(\a-2) + \s*\t );
							py(\a,\s,\t) = 0.5*( exp(\a*\t)/\a + \s*exp(2*\t)/2 - exp((\a-2)*\t)/(\a-2) - \s*\t );
						}
						]
						\tikzset{asymp/.style={dashed,gray!60,line width=0.4pt}}
						\begin{axis}[width=4.5cm,height=5.5cm,at={(0,0)},anchor=south west,
							axis lines=middle,ticks=none,xmin=-9,xmax=4,ymin=-1,ymax=6,clip=true]
							\addplot[blue,thick,domain=-1.2:1.2,samples=300,variable=\t]
							({px(-1,1,\t)}, {py(-1,1,\t)});
						\end{axis}
						\node[anchor=north, font=\small] at ([yshift=-0.4cm]current axis.south) {(i) $\gamma_+$, $a=-1$};
						\begin{axis}[width=4.5cm,height=5.5cm,at={(5cm,0)},anchor=south west,
							axis lines=middle,ticks=none,xmin=-13,xmax=9,ymin=-1,ymax=13,clip=true]
							\addplot[blue,thick,domain=-3:1.5,samples=300,variable=\t]
							({px(1,1,\t)}, {py(1,1,\t)});
							\draw[asymp] (axis cs:0,0) -- (axis cs:-11,11);
							\node[font=\tiny,gray] at (axis cs:-8.5,11.5) {$x+y=0$};
						\end{axis}
						\node[anchor=north, font=\small] at ([yshift=-0.4cm]current axis.south) {(ii) $\gamma_+$, $a=1$};
						\begin{axis}[width=4.5cm,height=5.5cm,at={(10cm,0)},anchor=south west,
							axis lines=middle,ticks=none,xmin=-2.5,xmax=8,ymin=-1,ymax=4,clip=true]
							\addplot[blue,thick,domain=-2.5:1.0,samples=300,variable=\t]
							({px(3,1,\t)}, {py(3,1,\t)});
							\draw[asymp] (axis cs:0,0) -- (axis cs:-2,2);
							\node[font=\tiny,gray] at (axis cs:-1.6,2.6) {$x+y=0$};
						\end{axis}
						\node[anchor=north, font=\small] at ([yshift=-0.4cm]current axis.south) {(iii) $\gamma_+$, $a=3$};
						\begin{axis}[width=4.5cm,height=5.5cm,at={(0,-6.3cm)},anchor=south west,
							axis lines=middle,ticks=none,xmin=-8,xmax=1,ymin=-3,ymax=5,clip=true]
							\addplot[blue,thick,domain=-1.2:1.2,samples=300,variable=\t]
							({px(-1,-1,\t)}, {py(-1,-1,\t)});
							\addplot[only marks,mark=*,mark size=1.5pt,red] coordinates {(-0.917,-0.583)};
						\end{axis}
						\node[anchor=north, font=\small] at ([yshift=-0.4cm]current axis.south) {(iv) $\gamma_-$, $a=-1$};
						\begin{axis}[width=4.5cm,height=5.5cm,at={(5cm,-6.3cm)},anchor=south west,
							axis lines=middle,ticks=none,xmin=-10,xmax=1.5,ymin=-3,ymax=10,clip=true]
							\addplot[blue,thick,domain=-3:1.5,samples=300,variable=\t]
							({px(1,-1,\t)}, {py(1,-1,\t)});
							\draw[asymp] (axis cs:0,0) -- (axis cs:-9,9);
							\node[font=\tiny,gray] at (axis cs:-7,9.3) {$x+y=0$};
							\addplot[only marks,mark=*,mark size=1.5pt,red] coordinates {(-0.25,0.75)};
						\end{axis}
						\node[anchor=north, font=\small] at ([yshift=-0.4cm]current axis.south) {(v) $\gamma_-$, $a=1$};
						\begin{axis}[width=4.5cm,height=5.5cm,at={(10cm,-6.3cm)},anchor=south west,
							axis lines=middle,ticks=none,xmin=-0.6,xmax=3,ymin=-2.2,ymax=1.2,clip=true]
							\addplot[blue,thick,domain=-2.5:1.0,samples=300,variable=\t]
							({px(3,-1,\t)}, {py(3,-1,\t)});
							\draw[asymp] (axis cs:0,0) -- (axis cs:2.2,-2.2);
							\node[font=\tiny,gray] at (axis cs:1.6,-2.0) {$x+y=0$};
							\addplot[only marks,mark=*,mark size=1.5pt,red] coordinates {(0.4167,-0.5833)};
						\end{axis}
						\node[anchor=north, font=\small] at ([yshift=-0.4cm]current axis.south) {(vi) $\gamma_-$, $a=3$};
					\end{tikzpicture}
					\caption{The mixed solitons $\gamma_+$ (top) and $\gamma_-$ (bottom) of \eqref{eq:gpm} for $a=b=-1,1,3$. Cusps are marked in red; the dashed line is the asymptote $x+y=0$.}
					\label{fig:mixedpm}
				\end{figure}
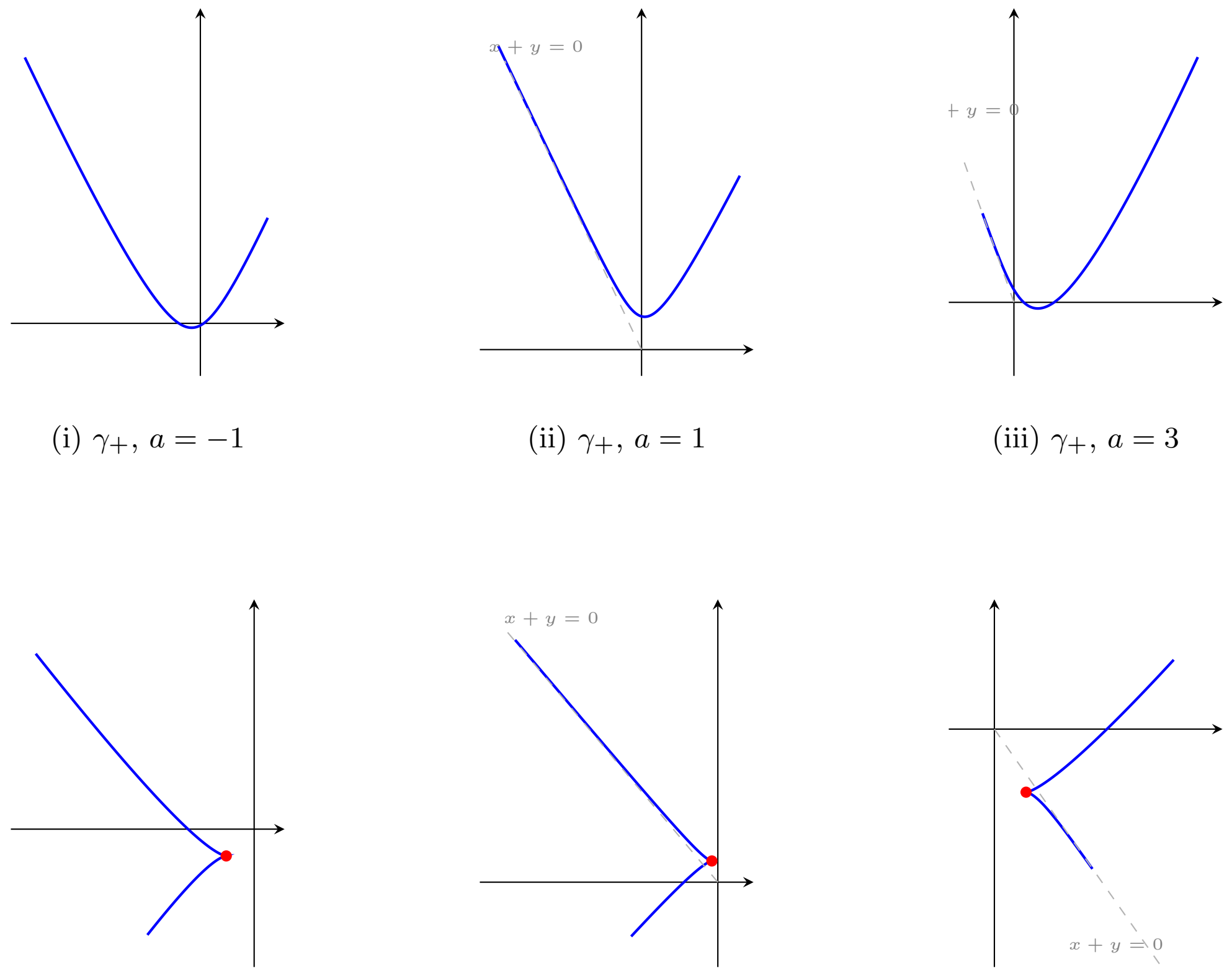

				\paragraph{The classes (I) and (II)}
				Normalizing $\frac{1}{k}=e^{(a-1)\theta}\pm e^{\theta}$ and integrating \eqref{eq:nullder} we obtain, up to a translation, the two curves
				\begin{equation}\label{eq:gpm}
					\gamma_{\pm}=\Bigl(\frac{e^{a\theta}}{a}\pm\frac{e^{2\theta}}{2}\Bigr)h_+
					+\Bigl(\frac{e^{(a-2)\theta}}{a-2}\pm\theta\Bigr)h_-,
				\end{equation}
				that is, after multiplying by $4a(a-2)$,
				\begin{align*}
					x&=2(a-2)e^{a\theta}\pm a(a-2)e^{2\theta}+2ae^{(a-2)\theta}\pm2a(a-2)\theta,\\
					y&=2(a-2)e^{a\theta}\pm a(a-2)e^{2\theta}-2ae^{(a-2)\theta}\mp2a(a-2)\theta.
				\end{align*}
				The arc length differential formula is $ds=\bigl|e^{(a-1)\theta}\pm e^{\theta}\bigr|d\theta$, and we read off the following.
				\begin{itemize}
					\item $\gamma_+$ is smooth for every $a$, whereas $\gamma_-$ has exactly one singular point, at $\theta=0$, where $e^{(a-2)\theta}=1$; there $\gamma_{\theta}$ vanishes and the trace has a cusp.
					\item The arc length is infinite as $\theta\to+\infty$ for every $a$, because the larger exponent $\max\{a-1,1\}$ is positive. As $\theta\to-\infty$ it is governed by the smaller exponent $\min\{a-1,1\}$ and diverges if and only if $a\le1$. Hence $\gamma_{\pm}$ is complete if and only if $a\le1$.
					\item If $a>0$, then $\hat{x}\to0$ and $|\hat{y}|\to\infty$ as $\theta\to-\infty$, so that $\gamma_{\pm}$ is asymptotic to the light-like line $x+y=0$. If $a<0$, both $\hat{x}$ and $\hat{y}$ diverge as $\theta\to-\infty$ and there is no asymptotic line.
				\end{itemize}
				In particular the curves $\gamma_+$ with $a\le1$, $a\ne0$, form a one-parameter family of complete smooth solutions, whose endpoint $a=1$ is the curve with $\frac{1}{k}=1+e^{\theta}$, namely $\gamma_+=(e^{\theta}+\frac{1}{2}e^{2\theta})h_++(\theta-e^{-\theta})h_-$; the curves $\gamma_-$ with $a\le1$ are complete but singular. Figure~\ref{fig:mixedpm} displays $\gamma_{\pm}$ for $a=-1,1,3$.
				
				\paragraph{The class (III)}
				If $\alpha=0$, that is $K=0$, then $\frac{1}{k}=\beta e^{\theta}$, $\hat{x}=\frac{C_-}{2(a-2)}e^{2\theta}$ and $\hat{y}=\frac{C_-}{a-2}\theta$ up to translations. This is, for every value of $a$, the curve \eqref{log sprial}, i.e.\ the light-like translating soliton $(\ref{trC(1,1)})$: smooth, but not complete. Thus the curve \eqref{log sprial} is simultaneously a translating soliton and a mixed soliton for every $a=b\ne0$, in the same way as the hyperbolic logarithmic spirals of Proposition \ref{hyper log spiral} satisfy \eqref{dr} for a one-parameter family of values of $(a,b)$.
				
				\paragraph{The class (IV)}
				If $\beta=0$, that is $C_-=0$, then $C$ is parallel to $(1,1)$ and, by Lemma \ref{lem:normalisation}(i), a translation removes it altogether: the equation reduces to \eqref{dr} with $b=a$, which is the hyperbolic case of Section \ref{sec5} on the exceptional line $b=a$ of Remark \ref{rem:nulltrans}. Indeed $\frac{1}{k}=\alpha e^{(a-1)\theta}$ is the hyperbolic logarithmic spiral $\gamma=(a-1-h)e^{(a-1+h)\theta}$ of Proposition \ref{hyper log spiral} with $l=a-1\ne\pm1$. It is smooth, and it is complete if and only if $a=1$, in which case $\frac{1}{k}$ is constant and the curve is the standard hyperbola $y^{2}-x^{2}=\alpha^{2}$.
				
\subsection{Summary}
From the preceding analysis, we can draw the following conclusions:

\begin{thm}[Classification of the self-similar solutions with $a^2 = b^2$]\label{thm:6.1}
	Consider the full self-similar equation \eqref{meq} with $a^2=b^2$ and $a\ne0$. By the reduction at the beginning of this section we may assume $a=b$, and by Lemma \ref{lem:normalisation} the translation vector $C$ enters only through $C_-=C_x-C_y$, i.e.\ through its class modulo the light-like direction $(1,1)$. Up to homotheties and isometries of $\mathbb{R}^{1,1}$ the solutions are the following.
	\begin{description}
		\item[\textbf{Case A ($a = b = 2$):}] exactly two solutions. For $C_-\ne0$ the curve \eqref{eq:a2normal}, which has one cusp and is asymptotic to the light-like line $x+y=0$; for $C_-=0$ the smooth curve \eqref{log sprial}. Both are incomplete, so no complete solution exists in this resonant case.
		\item[\textbf{Case B ($a = b \neq 0,2$):}] writing $\frac{1}{k}=\alpha e^{(a-1)\theta}+\beta e^{\theta}$ as in \eqref{eq:kab}, exactly four solutions for every value of $a$:
		\begin{itemize}
			\item[(I)] $\alpha\beta>0$: the curve $\gamma_+$ of \eqref{eq:gpm}, smooth for every $a$ and complete if and only if $a\le1$;
			\item[(II)] $\alpha\beta<0$: the curve $\gamma_-$ of \eqref{eq:gpm}, with exactly one cusp, complete if and only if $a\le1$;
			\item[(III)] $\alpha=0$: the light-like translating soliton \eqref{log sprial}, smooth and incomplete;
			\item[(IV)] $\beta=0$, i.e.\ $C_-=0$: the hyperbolic logarithmic spiral $(a-1-h)e^{(a-1+h)\theta}$, smooth, and complete if and only if $a=1$, when it is the standard hyperbola.
		\end{itemize}
	\end{description}
\end{thm}

\begin{cor}[Complete smooth mixed solitons]\label{cor:6.3}
	Up to homotheties and isometries, the solutions of \eqref{meq} with $a^{2}=b^{2}$, $a\ne0$, that are both complete and smooth are exactly
	\begin{enumerate}
		\item the curves $\gamma_+$ of \eqref{eq:gpm} with $a=b\le1$, $a\ne0$, and
		\item the standard hyperbola, arising for $a=b=1$ and $C_-=0$.
	\end{enumerate}
	In particular there is no complete smooth solution with $a=b>1$, and the solutions in the first item are not isolated: since $a$ cannot be normalized away (see the beginning of this section), they form a genuine one-parameter family.
\end{cor}

\section{Higher-dimensional consequences}\label{sec7}

The classification obtained in Sections \ref{sec4}--\ref{sec6} is not confined to the plane. In this section we show that it governs, in every dimension, two natural classes of self-similar solutions of the inverse mean curvature flow in Lorentz--Minkowski space: it classifies the translationally invariant ones \emph{completely} (Theorem \ref{thm:cyl}), and it reduces the rotationally invariant ones to an explicit ODE system (Proposition \ref{pro:rot}), whose distinguished solution is the homothetically expanding hyperbolic space $\mathbb{H}^n(\rho)$ (Example \ref{ex:hyperboloid}).

Throughout this section we set
\begin{equation}\label{eq:Rn1}
	\mathbb{R}^{n,1}:=\mathbb{R}^{1,1}\times\mathbb{R}^{n-1},
	\qquad
	\langle U,V\rangle:=\bigl(u_1v_1-u_2v_2\bigr)+\sum_{j=1}^{n-1}\tilde{u}_j\tilde{v}_j,
\end{equation}
for $U=(u_1,u_2,\tilde{u})$ and $V=(v_1,v_2,\tilde{v})$. Thus $\mathbb{R}^{n,1}$ is the $(n+1)$-dimensional Lorentz--Minkowski space, the second coordinate being the time-like one, and $\mathbb{R}^{1,1}$ is embedded as the coordinate plane $\tilde{u}=0$. For a non-degenerate hypersurface $\Sigma\subset\mathbb{R}^{n,1}$ with unit normal $\mathbf{n}$ we write $\mathbf{H}_\Sigma$ for its mean curvature vector, and we let $\mathcal{H}_\Sigma$ denote the sum of the principal curvatures of $\Sigma$ with respect to $\mathbf{n}$, normalized exactly as the curvature $k$ of Section \ref{sec3}, so that $\langle\mathbf{H}_\Sigma,\mathbf{n}\rangle=-\mathcal{H}_\Sigma$ and $\mathcal{H}_\Sigma=k$ when $n=1$. With this normalization the inverse mean curvature flow \eqref{eq:imcf} takes the form
\begin{equation}\label{eq:imcfn}
	\frac{\partial X}{\partial t}=\frac{1}{\mathcal{H}_\Sigma}\,\mathbf{n} ,
\end{equation}
which for $n=1$ is literally the flow equation of Section \ref{sec3}. Consequently every computation of that section applies verbatim to $\Sigma$ once $k$ is replaced by $\mathcal{H}_\Sigma$; in particular no sign convention beyond that of Section \ref{sec3} is used anywhere below.

\subsection{Cylinders over planar solitons}\label{subsec:cyl}

\begin{defi}\label{def:cyl}
	Let $\gamma:I\to\mathbb{R}^{1,1}$ be a regular curve which is everywhere space-like or everywhere time-like, with $\varepsilon=\langle\mathbf{T},\mathbf{T}\rangle_m=\pm1$ as in Section \ref{sec3}. Its \emph{cylinder} is the hypersurface
	\begin{equation*}
		\Sigma_\gamma:=\gamma(I)\times\mathbb{R}^{n-1}\subset\mathbb{R}^{n,1},
		\qquad
		X(s,\mathbf{u})=(\gamma(s),\mathbf{u}),\quad \mathbf{u}\in\mathbb{R}^{n-1},
	\end{equation*}
	where $s$ is the arc-length parameter of $\gamma$.
\end{defi}

\begin{lem}\label{lem:cylH}
	The cylinder $\Sigma_\gamma$ is a non-degenerate hypersurface with induced metric
	$g_{\Sigma}=\varepsilon\,ds^2+\sum_{j}du_j^2$ and unit normal
	$\mathbf{n}(s,\mathbf{u})=(\mathbf{N}(s),\mathbf{0})$, $\langle\mathbf{n},\mathbf{n}\rangle=-\varepsilon$. Its principal curvatures with respect to $\mathbf{n}$ are $k,0,\dots,0$, and
	\begin{equation}\label{eq:cylH}
		\mathbf{H}_{\Sigma_\gamma}=\mathbf{H}_\gamma,
		\qquad
		\mathcal{H}_{\Sigma_\gamma}=k .
	\end{equation}
	In particular $\Sigma_\gamma$ is space-like if and only if $\gamma$ is space-like.
\end{lem}

\begin{proof}
	The coordinate fields are $X_s=(\mathbf{T},\mathbf{0})$ and $X_{u_j}=(\mathbf{0},\mathbf{e}_j)$, whence $\langle X_s,X_s\rangle=\varepsilon$, $\langle X_{u_i},X_{u_j}\rangle=\delta_{ij}$ and $\langle X_s,X_{u_j}\rangle=0$, which gives $g_\Sigma$. Since $\langle\mathbf{N},\mathbf{T}\rangle_m=0$ and $\langle\mathbf{N},\mathbf{N}\rangle_m=-\varepsilon$ by $\mathbf{N}=h\mathbf{T}$, the field $\mathbf{n}=(\mathbf{N},\mathbf{0})$ is a unit normal.
	
	The ambient connection is the flat one. The Frenet equation $\mathbf{T}_s=k\mathbf{N}$ gives $\nabla_{X_s}X_s=(k\mathbf{N},\mathbf{0})=k\,\mathbf{n}$, while $\nabla_{X_{u_i}}X_{u_j}=0$ and $\nabla_{X_s}X_{u_j}=0$, so the $\mathbb{R}^{n-1}$ directions are totally geodesic and contribute nothing to the trace. Taking the trace of the second fundamental form with respect to $g_\Sigma$ therefore yields $\mathbf{H}_{\Sigma_\gamma}=\mathbf{H}_\gamma$, and $\mathcal{H}_{\Sigma_\gamma}=k$, the curvature of $\gamma$. In particular no term of the type $(n-1)/r$ appears in $\mathcal{H}_{\Sigma_\gamma}$, in contrast with the rotational hypersurfaces of Subsection \ref{subsec:rot}; this is the reason why the planar theory transplants verbatim.
\end{proof}

\begin{thm}\label{thm:cyl}
	Let $n\ge 2$ and let $\gamma:I\to\mathbb{R}^{1,1}$ be as in Definition \ref{def:cyl}.
	\begin{enumerate}
		\item[\rm(i)] $\gamma$ satisfies \eqref{meq} with data $(a,b,C)$ if and only if $\Sigma_\gamma$ generates a self-similar solution of the inverse mean curvature flow in $\mathbb{R}^{n,1}$, that is, $\Sigma_t=\Phi_t(\Sigma_\gamma)$ solves \eqref{eq:imcfn} for some $C^1$ family of similarities $\Phi_t$ of $\mathbb{R}^{n,1}$ with $\Phi_0=\mathrm{id}$. The family may then be taken to be
		\begin{equation}\label{eq:Phit}
			\Phi_t(x,\mathbf{u}):=\bigl(g(t)e^{hf(t)}x+\mathbf{G}(t),\ g(t)\mathbf{u}\bigr),
		\end{equation}
		where $f,g,\mathbf{G}$ are as in Section \ref{sec3}.
		\item[\rm(ii)] Conversely, let $\Sigma\subset\mathbb{R}^{n,1}$ be a connected, properly embedded, non-degenerate hypersurface invariant under the translation group $\{(0,0,\mathbf{u}):\mathbf{u}\in\mathbb{R}^{n-1}\}\cong\mathbb{R}^{n-1}$. Then $\Sigma=\Sigma_\gamma$ for a curve $\gamma\subset\mathbb{R}^{1,1}$, unique up to reparametrization.
		\item[\rm(iii)] $\Sigma_\gamma$ is complete, resp.\ smooth, if and only if $\gamma$ is.
	\end{enumerate}
	Consequently, Theorem \ref{thm:main} classifies the $\mathbb{R}^{n-1}$-invariant self-similar solutions of the inverse mean curvature flow in $\mathbb{R}^{n,1}$ for every $n\ge2$.
\end{thm}

\begin{proof}
	(i) By Lemma \ref{lem:cylH}, $\mathbf{n}=(\mathbf{N},\mathbf{0})$ and $\mathcal{H}_{\Sigma_\gamma}=k$, so that \eqref{eq:imcfn} for the cylinders $\Sigma_{\hat{\gamma}(\cdot,t)}$, $\hat{\gamma}(\cdot,t)=g(t)e^{hf(t)}\gamma+\mathbf{G}(t)$, is the planar flow equation of Section \ref{sec3} for $\hat{\gamma}$; hence $\gamma$ solves \eqref{meq} if and only if $t\mapsto\Sigma_{\hat{\gamma}(\cdot,t)}$ solves the flow. The map \eqref{eq:Phit} is a similarity of $\mathbb{R}^{n,1}$ with $\Phi_t(\Sigma_\gamma)=\hat{\gamma}(I,t)\times g(t)\mathbb{R}^{n-1}=\Sigma_{\hat{\gamma}(\cdot,t)}$, and the parametrizations $\Phi_t\circ X$ and $(s,\mathbf{u})\mapsto(\hat{\gamma}(s,t),\mathbf{u})$ differ by the tangential diffeomorphism $\mathbf{u}\mapsto g(t)\mathbf{u}$. Thus a solution of \eqref{meq} generates the self-similar solution \eqref{eq:Phit}.

	Conversely, let $\Sigma_t=\Phi_t(\Sigma_\gamma)$ solve \eqref{eq:imcfn} for a $C^1$ family of similarities with $\Phi_0=\mathrm{id}$. Its velocity at $t=0$ is an infinitesimal similarity $V(X)=\mu X+LX+W$ with $\mu\in\mathbb{R}$, $L\in\mathfrak{so}(n,1)$, $W\in\mathbb{R}^{n,1}$, and as in the derivation of \eqref{cfeq} the flow equation at $t=0$ reads $\langle V(X),\mathbf{n}\rangle=-1/k$ at every $X=(\gamma(s),\mathbf{u})\in\Sigma_\gamma$. Decompose $L$ into blocks with respect to $\mathbb{R}^{1,1}\times\mathbb{R}^{n-1}$ and write $W=(W_1,W_2)$; then $L_{11}=a\,h\in\mathfrak{so}(1,1)$ for some $a\in\mathbb{R}$, and since $\mathbf{n}=(\mathbf{N},\mathbf{0})$ and $\langle h\gamma,\mathbf{N}\rangle_m=-\langle\gamma,\mathbf{T}\rangle_m$, the equation becomes
	\begin{equation*}
		a\langle\gamma,\mathbf{T}\rangle_m-\mu\langle\gamma,\mathbf{N}\rangle_m-\langle L_{12}\mathbf{u},\mathbf{N}\rangle_m-\langle W_1,\mathbf{N}\rangle_m=\frac{1}{k}.
	\end{equation*}
	The right-hand side does not depend on $\mathbf{u}$, so the term linear in $\mathbf{u}$ vanishes identically; as $k\ne0$, the normals $\mathbf{N}(s)$ span $\mathbb{R}^{1,1}$ and hence $L_{12}=0$. What remains is \eqref{meq} with $b=\mu$ and $C=W_1$, and $\Phi_t$ may be replaced by \eqref{eq:Phit}.

	(ii) Invariance under $\mathbb{R}^{n-1}$ means that $\Sigma$ is foliated by the affine planes $p+\{0\}\times\mathbb{R}^{n-1}$, so that the projection $\pi:\mathbb{R}^{n,1}\to\mathbb{R}^{1,1}$ has constant rank $1$ on $\Sigma$. As $\Sigma$ is connected and properly embedded, $\pi(\Sigma)$ is the trace of a regular curve $\gamma$, unique up to reparametrization, and $\Sigma=\pi(\Sigma)\times\mathbb{R}^{n-1}$. Non-degeneracy of $\Sigma$ gives $\langle\gamma',\gamma'\rangle_m\ne0$, hence of constant sign on the connected interval $I$, as required in Definition \ref{def:cyl}.

	(iii) The induced metric is the product metric of Lemma \ref{lem:cylH}; a product is complete if and only if each factor is, and $X$ is a smooth immersion exactly where $\gamma$ is.
\end{proof}

Time-like profile curves are covered as well: by Section \ref{sec3} they are the images of the space-like ones under the reflection $z\mapsto hz$, which for $n\ge2$ is not induced by an isometry of $\mathbb{R}^{n,1}$, so that the time-like cylinders $\Sigma_{h\gamma}$ are not congruent to the space-like ones.

\begin{rem}[Lorentzian grim-reaper cylinders]\label{rem:grim}
	Taking for $\gamma$ the complete smooth translating soliton of Theorem \ref{translating thm}, Theorem \ref{thm:cyl} produces, in each dimension, a complete translating soliton of the inverse mean curvature flow in $\mathbb{R}^{n,1}$. These are the exact analogues of the grim-reaper cylinders of the Euclidean mean curvature flow, which arise as blow-up limits at type-II singularities. In the same way, Theorem \ref{thm:cyl} together with Theorem \ref{thm:main} shows that any self-similar blow-up limit of the inverse mean curvature flow in $\mathbb{R}^{n,1}$ which is invariant under $\mathbb{R}^{n-1}$ must be one of the cylinders over the curves classified in Sections \ref{sec4}--\ref{sec6}.
\end{rem}

\subsection{Rotational hypersurfaces}\label{subsec:rot}

Cylinders are the hypersurfaces invariant under a group of translations. We now consider the other natural invariant class, that of hypersurfaces invariant under the rotation group about a time-like axis. For this subsection it is convenient to reorder the coordinates of \eqref{eq:Rn1} and write
\begin{equation*}
	\mathbb{R}^{n,1}=\bigl\{(x,y):x\in\mathbb{R}^n,\ y\in\mathbb{R}\bigr\},
	\qquad
	\bigl\langle(x,y),(x',y')\bigr\rangle=\langle x,x'\rangle_{\mathbb{R}^n}-yy',
\end{equation*}
so that the $y$-axis is time-like, the group $SO(n)$ acts on the $x$-variable, and the plane $\mathbb{R}^{1,1}$ is recovered as $\{(x_1,0,\dots,0,y)\}$ with $h$ acting by $(x_1,y)\mapsto(y,x_1)$.

\begin{defi}\label{def:rot}
	Let $\gamma=(r,y):I\to\mathbb{R}^{1,1}$ be a space-like unit-speed curve with $r>0$. The associated \emph{rotational hypersurface} is
	\begin{equation*}
		\Sigma^{\mathrm{rot}}_\gamma:=\bigl\{\bigl(r(s)\,\omega,\ y(s)\bigr):s\in I,\ \omega\in S^{n-1}\bigr\}\subset\mathbb{R}^{n,1}.
	\end{equation*}
\end{defi}

\begin{lem}\label{lem:rotH}
	Write $\mathbf{T}=(r',y')=e^{h\theta}$ and $\mathbf{N}=h\mathbf{T}=(y',r')$ as in Section \ref{sec3}. Then $\Sigma^{\mathrm{rot}}_\gamma$ is a space-like hypersurface with time-like unit normal $\mathbf{n}(s,\omega)=(y'\omega,r')$, and its principal curvatures with respect to $\mathbf{n}$ are
	\begin{equation}\label{eq:rotpc}
		\lambda_0=k,\qquad \lambda_1=\dots=\lambda_{n-1}=\frac{y'}{r},
	\end{equation}
	so that
	\begin{equation}\label{eq:rotH}
		\mathcal{H}_\Sigma=k+(n-1)\,\frac{y'}{r} .
	\end{equation}
\end{lem}

\begin{proof}
	One has $\langle\mathbf{n},\mathbf{n}\rangle=y'^2-r'^2=-1$ since $\gamma$ is space-like and unit speed, and $\mathbf{n}$ is orthogonal to $X_s=(r'\omega,y')$ and to the fields tangent to the spheres, so $\mathbf{n}$ is a time-like unit normal and $\Sigma^{\mathrm{rot}}_\gamma$ is space-like.
	
	The Frenet equations give $r''=ky'$ and $y''=kr'$, hence
	\begin{equation*}
		\partial_s\mathbf{n}=(y''\omega,r'')=(kr'\omega,ky')=k\,X_s ,
	\end{equation*}
	which is $\lambda_0=k$. For the spherical directions, let $v\in T_\omega S^{n-1}$ be a unit vector; the corresponding tangent vector of $\Sigma^{\mathrm{rot}}_\gamma$ is $(rv,0)$, of length $r$, while the derivative of $\mathbf{n}$ in that direction is $(y'v,0)$. Dividing by $r$ gives $\lambda_i=y'/r$ for each of the $n-1$ spherical directions. Summing yields \eqref{eq:rotH}.
\end{proof}

\begin{pro}\label{pro:rot}
	Let $(\Sigma_t)$ be a self-similar solution of the inverse mean curvature flow all of whose members are rotational hypersurfaces about one and the same time-like axis. Then the associated similarities preserve that axis, so that $a=0$ and $C=c\,e_y$ in the notation of Section \ref{sec3}. The profile curve $\gamma=(r,y)$ then satisfies
	\begin{equation}\label{eq:rotmeq}
		-b\,\langle \gamma,\mathbf{N}\rangle_m-\langle C,\mathbf{N}\rangle_m
		=\frac{1}{\,k+(n-1)\dfrac{y'}{r}\,},
	\end{equation}
	which is equation \eqref{meq} with $a=0$ and with the curvature $k$ replaced by $\mathcal{H}_\Sigma$. Equivalently, in terms of the hyperbolic tangent angle $\theta$,
	\begin{equation}\label{eq:rotode}
		r'=\cosh\theta,\qquad y'=\sinh\theta,\qquad
		\theta'=\frac{1}{b\,(y\cosh\theta-r\sinh\theta)+c\cosh\theta}-(n-1)\frac{\sinh\theta}{r}.
	\end{equation}
\end{pro}

\begin{proof}
	By hypothesis every $\Sigma_t$ is invariant under the same subgroup $SO(n)\subset O(n,1)$, so the similarity $\Phi_t$ carrying $\Sigma_0$ to $\Sigma_t$ normalizes $SO(n)$ and therefore fixes its axis, the $y$-axis, as a set. A non-trivial boost has no time-like fixed direction and hence cannot do this, so $\Phi_t$ is a composition of a homothety, a translation along the $y$-axis and an element of $SO(n)$. Thus $f\equiv0$, i.e.\ $a=0$, and $\mathbf{G}(t)$ is parallel to $e_y$. We emphasize that the invariance of \emph{all} the $\Sigma_t$ about one fixed axis is a hypothesis, not a normalization: a totally umbilic solution such as the one in Example \ref{ex:hyperboloid} is invariant under every rotation subgroup of $O(n,1)$, so its axis is not determined by the hypersurface.
	
	Writing $\hat{X}(t)=g(t)X+\mathbf{G}(t)$ and differentiating at $t=0$ gives $\partial_t\hat{X}|_{0}=bX+c\,e_y$ with $b=g'(0)$, $c=\mathbf{G}'(0)\cdot e_y$. We now repeat the computation of Section \ref{sec3} verbatim, with $k$ replaced by $\mathcal{H}_\Sigma$; by Lemma \ref{lem:rotH},
	\begin{equation*}
		\langle X,\mathbf{n}\rangle=\bigl\langle(r\omega,y),(y'\omega,r')\bigr\rangle=ry'-yr'=\langle\gamma,\mathbf{N}\rangle_m,
		\qquad
		\langle e_y,\mathbf{n}\rangle=-r' ,
	\end{equation*}
	so that $\langle C,\mathbf{N}\rangle_m=-cr'=\langle c\,e_y,\mathbf{n}\rangle$. The pairings entering the derivation of \eqref{meq} are therefore exactly those of Section \ref{sec3}, and one obtains \eqref{eq:rotmeq}, that is, equation \eqref{meq} with $a=0$ and with $k$ replaced by $\mathcal{H}_\Sigma$. Solving \eqref{eq:rotmeq} for $k=\theta'$ and using $\mathbf{T}=e^{h\theta}$ yields \eqref{eq:rotode}.
\end{proof}

\begin{example}\label{ex:hyperboloid}
	Let $\gamma(s)=\rho\bigl(\sinh\frac{s}{\rho},\cosh\frac{s}{\rho}\bigr)$, $s>0$, so that $\Sigma^{\mathrm{rot}}_\gamma$ is the hyperbolic space
	\begin{equation*}
		\mathbb{H}^n(\rho)=\bigl\{(x,y)\in\mathbb{R}^{n,1}:|x|^2-y^2=-\rho^2,\ y>0\bigr\},
	\end{equation*}
	a complete totally umbilic space-like hypersurface. (Strictly speaking, Definition \ref{def:rot} requires $r>0$, so $\Sigma^{\mathrm{rot}}_\gamma$ is $\mathbb{H}^n(\rho)$ with the vertex $x=0$, $y=\rho$ removed; as always for a rotational parametrization, the hypersurface extends smoothly across the axis, and we identify the two.) Here $k=1/\rho$, $\langle\gamma,\mathbf{N}\rangle_m=-\rho$ and $y'/r=1/\rho$, so $\mathcal{H}_\Sigma=n/\rho$ by \eqref{eq:rotH}. With $c=0$, equation \eqref{eq:rotmeq} becomes $b\rho=\rho/n$, that is
	\begin{equation}\label{eq:hyperboloidrate}
		b=\frac1n,\qquad \rho(t)=\rho_0\,e^{t/n} .
	\end{equation}
	Thus $\mathbb{H}^n(\rho)$ is a homothetically expanding soliton, the exact Lorentzian counterpart of the expanding round sphere of the Euclidean inverse mean curvature flow. For $n=1$ formula \eqref{eq:hyperboloidrate} returns $b=1$, which is the hyperbola of the planar classification; this provides an independent check of the conventions of Section \ref{sec3}.
\end{example}


\section*{Declarations}
\begin{itemize}
	\item Funding: This work was supported by the Sichuan Provincial Natural Science Foundation Youth Science Fund Project, Category B (Grant No.\ 2026NSFSC0707).
	\item Conflict of interest: The author declares no conflict of interest.
	\item Ethics approval and consent to participate: Not applicable.
	\item Consent for publication: Not applicable.
	\item Data availability: This manuscript has no associated data.
	\item Materials availability: Not applicable.
	\item Code availability: Not applicable.
	\item Author contribution: The author is solely responsible for the manuscript.
\end{itemize}

\bibliography{bibouli}


\begin{thebibliography}{13}
\ifx \bisbn   \undefined \def \bisbn  #1{ISBN #1}\fi
\ifx \binits  \undefined \def \binits#1{#1}\fi
\ifx \bauthor  \undefined \def \bauthor#1{#1}\fi
\ifx \batitle  \undefined \def \batitle#1{#1}\fi
\ifx \bjtitle  \undefined \def \bjtitle#1{#1}\fi
\ifx \bvolume  \undefined \def \bvolume#1{\textbf{#1}}\fi
\ifx \byear  \undefined \def \byear#1{#1}\fi
\ifx \bissue  \undefined \def \bissue#1{#1}\fi
\ifx \bfpage  \undefined \def \bfpage#1{#1}\fi
\ifx \blpage  \undefined \def \blpage #1{#1}\fi
\ifx \burl  \undefined \def \burl#1{\textsf{#1}}\fi
\ifx \doiurl  \undefined \def \doiurl#1{\url{https://doi.org/#1}}\fi
\ifx \betal  \undefined \def \betal{\textit{et al.}}\fi
\ifx \binstitute  \undefined \def \binstitute#1{#1}\fi
\ifx \binstitutionaled  \undefined \def \binstitutionaled#1{#1}\fi
\ifx \bctitle  \undefined \def \bctitle#1{#1}\fi
\ifx \beditor  \undefined \def \beditor#1{#1}\fi
\ifx \bpublisher  \undefined \def \bpublisher#1{#1}\fi
\ifx \bbtitle  \undefined \def \bbtitle#1{#1}\fi
\ifx \bedition  \undefined \def \bedition#1{#1}\fi
\ifx \bseriesno  \undefined \def \bseriesno#1{#1}\fi
\ifx \blocation  \undefined \def \blocation#1{#1}\fi
\ifx \bsertitle  \undefined \def \bsertitle#1{#1}\fi
\ifx \bsnm \undefined \def \bsnm#1{#1}\fi
\ifx \bsuffix \undefined \def \bsuffix#1{#1}\fi
\ifx \bparticle \undefined \def \bparticle#1{#1}\fi
\ifx \barticle \undefined \def \barticle#1{#1}\fi
\bibcommenthead
\ifx \bconfdate \undefined \def \bconfdate #1{#1}\fi
\ifx \botherref \undefined \def \botherref #1{#1}\fi
\ifx \url \undefined \def \url#1{\textsf{#1}}\fi
\ifx \bchapter \undefined \def \bchapter#1{#1}\fi
\ifx \bbook \undefined \def \bbook#1{#1}\fi
\ifx \bcomment \undefined \def \bcomment#1{#1}\fi
\ifx \oauthor \undefined \def \oauthor#1{#1}\fi
\ifx \citeauthoryear \undefined \def \citeauthoryear#1{#1}\fi
\ifx \endbibitem  \undefined \def \endbibitem {}\fi
\ifx \bconflocation  \undefined \def \bconflocation#1{#1}\fi
\ifx \arxivurl  \undefined \def \arxivurl#1{\textsf{#1}}\fi
\csname PreBibitemsHook\endcsname

\bibitem[\protect\citeauthoryear{Urbas}{1991}]{MR1085136}
\begin{barticle}
\bauthor{\bsnm{Urbas}, \binits{J.I.E.}}:
\batitle{An expansion of convex hypersurfaces}.
\bjtitle{J. Differential Geom.}
\bvolume{33}(\bissue{1}),
\bfpage{91}--\blpage{125}
(\byear{1991})
\doiurl{10.4310/jdg/1214446031}
\end{barticle}
\endbibitem

\bibitem[\protect\citeauthoryear{Urbas}{1990}]{MR1082861}
\begin{barticle}
\bauthor{\bsnm{Urbas}, \binits{J.I.E.}}:
\batitle{On the expansion of starshaped hypersurfaces by symmetric functions of
  their principal curvatures}.
\bjtitle{Math. Z.}
\bvolume{205}(\bissue{3}),
\bfpage{355}--\blpage{372}
(\byear{1990})
\doiurl{10.1007/BF02571249}
\end{barticle}
\endbibitem

\bibitem[\protect\citeauthoryear{Gerhardt}{1990}]{Gerhardt1990}
\begin{barticle}
\bauthor{\bsnm{Gerhardt}, \binits{C.}}:
\batitle{Flow of nonconvex hypersurfaces into spheres}.
\bjtitle{J. Differential Geom.}
\bvolume{32}(\bissue{1}),
\bfpage{299}--\blpage{314}
(\byear{1990})
\doiurl{10.4310/jdg/1214445048}
\end{barticle}
\endbibitem

\bibitem[\protect\citeauthoryear{Huisken and Ilmanen}{2001}]{Huisken2001}
\begin{barticle}
\bauthor{\bsnm{Huisken}, \binits{G.}},
\bauthor{\bsnm{Ilmanen}, \binits{T.}}:
\batitle{The inverse mean curvature flow and the {R}iemannian {P}enrose
  inequality}.
\bjtitle{J. Differential Geom.}
\bvolume{59}(\bissue{3}),
\bfpage{353}--\blpage{437}
(\byear{2001})
\doiurl{10.4310/jdg/1090349447}
\end{barticle}
\endbibitem

\bibitem[\protect\citeauthoryear{Huisken and
  Ilmanen}{1997}]{huisken1997riemannian}
\begin{botherref}
\oauthor{\bsnm{Huisken}, \binits{G.}},
\oauthor{\bsnm{Ilmanen}, \binits{T.}}:
The riemannian penrose inequality.
Internat. Math. Res. Notices.
\textbf{1997}(20)
(1997)
\doiurl{10.1155/S1073792897000664}
\end{botherref}
\endbibitem

\bibitem[\protect\citeauthoryear{Bray and Neves}{2004}]{MR2081443}
\begin{barticle}
\bauthor{\bsnm{Bray}, \binits{H.L.}},
\bauthor{\bsnm{Neves}, \binits{A.}}:
\batitle{Classification of prime 3-manifolds with {Y}amabe invariant greater
  than {$\Bbb{RP}^3$}}.
\bjtitle{Ann. of Math. (2)}
\bvolume{159}(\bissue{2}),
\bfpage{887}
(\byear{2004})
\doiurl{10.4007/annals.2004.159.887}
\end{barticle}
\endbibitem

\bibitem[\protect\citeauthoryear{Drugan et~al.}{2016}]{Drugan2016}
\begin{barticle}
\bauthor{\bsnm{Drugan}, \binits{G.}},
\bauthor{\bsnm{Lee}, \binits{H.}},
\bauthor{\bsnm{Wheeler}, \binits{G.}}:
\batitle{Solitons for the inverse mean curvature flow}.
\bjtitle{Pac. J. Math.}
\bvolume{284}(\bissue{2}),
\bfpage{309}--\blpage{326}
(\byear{2016})
\doiurl{10.2140/pjm.2016.284.309}
\end{barticle}
\endbibitem

\bibitem[\protect\citeauthoryear{Urbas}{1999}]{MR1793674}
\begin{barticle}
\bauthor{\bsnm{Urbas}, \binits{J.}}:
\batitle{Convex curves moving homothetically by negative powers of their
  curvature}.
\bjtitle{Asian J. Math.}
\bvolume{3}(\bissue{3}),
\bfpage{635}--\blpage{656}
(\byear{1999})
\doiurl{10.4310/AJM.1999.v3.n3.a4}
\end{barticle}
\endbibitem

\bibitem[\protect\citeauthoryear{Chang}{2020}]{Chang2020}
\begin{barticle}
\bauthor{\bsnm{Chang}, \binits{J.-E.}}:
\batitle{Self-similar solutions to the inverse mean curvature flow in
  $\mathbb{R}^2$}.
\bjtitle{Results Math.}
\bvolume{75}(\bissue{4}),
\bfpage{175}--\blpage{19}
(\byear{2020})
\doiurl{10.1007/s00025-020-01304-2}
\end{barticle}
\endbibitem

\bibitem[\protect\citeauthoryear{Halldorsson}{2015}]{Halldorsson2015}
\begin{botherref}
\oauthor{\bsnm{Halldorsson}, \binits{H.P.}}:
Self-similar solutions to the mean curvature flow in the {M}inkowski plane
  {$\mathbb{R}^{1,1}$}.
J. Reine Angew. Math.,
209--243
(2015)
\doiurl{10.1515/crelle-2013-0054}
\end{botherref}
\endbibitem

\bibitem[\protect\citeauthoryear{da~Silva and Tenenblat}{2023}]{daSilva2023}
\begin{barticle}
\bauthor{\bsnm{Silva}, \binits{F.N.}},
\bauthor{\bsnm{Tenenblat}, \binits{K.}}:
\batitle{Self-similar solutions to the curvature flow and its inverse on the
  2-dimensional light cone}.
\bjtitle{Ann. Mat. Pura Appl.}
\bvolume{202}(\bissue{1}),
\bfpage{253}--\blpage{285}
(\byear{2023})
\doiurl{10.1007/s10231-022-01240-8}
\end{barticle}
\endbibitem

\bibitem[\protect\citeauthoryear{Sobczyk}{1995}]{hypernumber}
\begin{barticle}
\bauthor{\bsnm{Sobczyk}, \binits{G.}}:
\batitle{The hyperbolic number plane}.
\bjtitle{COLLEGE MATH. J.}
\bvolume{26}(\bissue{4}),
\bfpage{268}--\blpage{280}
(\byear{1995})
\doiurl{10.1080/07468342.1995.11973712}
\end{barticle}
\endbibitem

\bibitem[\protect\citeauthoryear{Catoni et~al.}{2011}]{MR2797703}
\begin{bbook}
\bauthor{\bsnm{Catoni}, \binits{F.}},
\bauthor{\bsnm{Boccaletti}, \binits{D.}},
\bauthor{\bsnm{Cannata}, \binits{R.}},
\bauthor{\bsnm{Catoni}, \binits{V.}},
\bauthor{\bsnm{Zampetti}, \binits{P.}}:
\bbtitle{Geometry of {M}inkowski Space-time}.
\bsertitle{SpringerBriefs in Physics},
p. \bfpage{116}.
\bpublisher{Springer}, \blocation{???}
(\byear{2011}).
\doiurl{10.1007/978-3-642-17977-8} .
\burl{https://doi.org/10.1007/978-3-642-17977-8}
\end{bbook}
\endbibitem

\end{thebibliography}
\end{document}